\documentclass[journal, onecolumn, 11pt, draftclsnofoot]{IEEEtran}

\ifCLASSINFOpdf
\else
\fi

\usepackage{amsmath,amssymb,amsfonts}
\usepackage{cite}
\usepackage{url}
\usepackage{xcolor}

\usepackage{nicefrac}
\usepackage{mathtools}
\usepackage{hyperref}
\usepackage{enumitem}
\usepackage[pdftex]{graphicx}

\begin{document}
%
\title{Solving non-linear Kolmogorov equations in large dimensions by using deep learning: a numerical comparison of discretization schemes}
%
%
%



\author{ \href{https://orcid.org/0000-0002-2311-4380}{\includegraphics[scale=0.06]{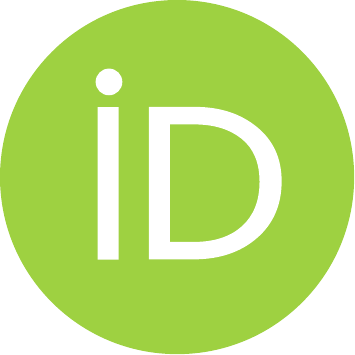}\hspace{1mm} Raffaele Marino} \\
	Dipartimento di Fisica,\\
	Sapienza Università di Roma\\
	Piazzale Aldo Moro 5, 00185, Rome, Italy \\
	\texttt{raffaele.marino@uniroma1.it} \\
	
	Nicolas Macris \\
	School of Computer and Communication Science,\\
	\'{E}cole Polytechnique F\'{e}d\'{e}rale de Lausanne.\\
	Rte Cantonale, 1015 Lausanne, Switzerland\\

}

\maketitle

\begin{abstract}
Non-linear partial differential Kolmogorov equations are successfully used to describe a wide range of time dependent phenomena, in natural sciences, engineering or even finance. For example, in physical systems, the Allen-Cahn equation describes pattern formation associated to phase transitions. In finance, instead, the Black-Scholes equation describes the evolution of the price of derivative investment instruments. Such modern applications often require to solve these equations in high-dimensional regimes in which classical approaches are ineffective. Recently, an interesting new approach based on deep learning has been introduced by E, Han, and Jentzen \cite{han2018solving, E2017}. The main idea is to construct a deep network which is trained from the samples of discrete stochastic differential equations underlying Kolmogorov's equation. The network is able to approximate, numerically at least, the solutions of the Kolmogorov equation with polynomial complexity in whole spatial domains.

In this contribution we study variants of the deep networks by using different discretizations schemes of the stochastic differential equation. We compare the performance of the associated networks, on benchmarked examples, and show that, for some discretization schemes, improvements in the accuracy are possible without affecting the observed computational complexity.
\end{abstract}


\section{Introduction}
\label{sec::introduc}
Algorithms based on the theory of deep learning have become essential in a wide variety scientific disciplines. 
In this paper we are concerned with applications to the solution of high-dimensional semi-linear parabolic partial differential equations (PDEs). The importance of such PDEs in finance, mathematics, natural science and engineering, cannot be understated and vast amounts of efforts have been deployed to develop numerical solution methods. However, in very high dimensions of the order of one-hundred or more many difficulties remain due to the curse of dimensionality \cite{zbMATH03126094}. This is the case for typical grid or mesh based methods, but also for Monte Carlo based approaches which only avoid the curse of dimensionality when computing the solution at (essentially) a single spatial point, but {\it not} for a whole domain. We refer to \cite{beck2019deep, DBLP:journals/jscic/Chan-Wai-NamMW19} for good up-to-date summaries of the state of the art as well as numerous references to old and modern methods. 


Using neural networks to solve and learn ODE's and PDE's is an old idea, see for example \cite{lee1990neural, meade1994solution, dissanayake1994neural, 712178, articleLagaris, MALEK2006260}. Recently however new successful ideas using the capabilities of deep learning have revived this approach in various directions. One such direction explores "physics informed" deep learning architectures that use a small or moderate amount of data and leverage on prior knowledge of the class of PDEs governing the physical dynamics, see for example \cite{raissi2019physics, 10.3389/fphy.2020.00042, articleScienceFluid}. On the front of very high-dimensional problems, 
 a mesh free "Deep Galerkin Method" \cite{sirignano2018dgm} based on deep neural networks trained with quadratic losses constructed out of the differential operator, initial and boundary conditions, is shown to be successful for spatial dimensions of the order of $d=100$ - $200$. Still for very high dimensional problems, in another direction which is our interest here, E, Han and Jentzen \cite{han2018solving, E2017} proposed a very interesting methodology, based on the connection between non-linear backward Kolmogorov equations and forward-backward stochastic differential equations (FBSDE). The general link between FBSDE, non-linear backward Kolmogorov equations, and the non-linear Feynman-Kac formula, goes back to the pioneering theory of Bismut \cite{zbMATH03506547} and Pardoux and Peng \cite{PARDOUX199055} on stochastic optimal control. The backward equation is used to construct neural network architectures trained from forward stochastic trajectories so as to minimize a quadratic error with respect to the terminal condition of the PDE problem. It has been possible to approximate the solution - at a single time and whole space domain - for spatial dimensions of the orders of $d=100$ 
with polynomial complexity. The approximations are benchmarked for many PDE's for which there are trusted numerical or sometimes exact solutions (see \cite{zhou2021actor} for an up to date summary). This approach has been further explored in various other directions, e.g.,  better architectures and loss functions, convergence, handling of memory problems, more severe non-linearities, approximating the solution at all times with a single set of network weights, even higher dimensions, partial integrodifferential equations \cite{BeckStillanother2019, hure2019some, DBLP:journals/corr/abs-1908-00412, DBLP:journals/jscic/Chan-Wai-NamMW19, raissi2018forward, beck2019deep, GS20-925, zhou2021actor, han2020convergence, jiang2021convergence}.  For an up to date summary of using deep learning for solving high dimensional PDEs one can see \cite{weinan2021algorithms} and reference therein.

Rigorous analysis that the deep learning methodology overcomes the curse of dimensionality in whole spatial domains is an open problem and rigorous results are still scarce. The representation power of deep learning architectures has been addressed in \cite{berner2018analysis, grohs2018proof} for linear equations, in  
\cite{Hutzenthaler-2020} for a class of semi-linear equations, and the convergence of the deep learning methodology for solving FBSDE has been adressed in \cite{han2020convergence, jiang2021convergence}. We also note that approximation algorithms which provably overcome the curse of dimensionality for non linear equations for single space-time points (and any time horizon) are based on multilevel Picard approximations \cite{Weinan2016MultilevelPI,
10.1007/s10915-018-00903-0} (see \cite{becker2020numerical} and references therein for the state of the art).


In this article we further investigate the original algorithm proposed in \cite{han2018solving, E2017} in a new direction which has not been investigated. To construct and train the neural network architecture it is necessary to discretize the underlying SDE's. Till now only the simplest Euler-Maruyama discretization has been used. We explore different discretizations, beyond the simplest Euler-Maruyama scheme. In particular, we show that a discretization of higher order, such as a multidimensional version of Milstein's second-order scheme \cite{kloeden2013numerical}, can markedly improve the accuracy of the algorithm without affecting the convergence time and computational complexity.  We also study a first-order scheme put forward by Leimkuhler and Matthews \cite{leimkuhler2012rational} in the context of Langevin dynamics and sampling, and show that its practical application to the equations considered here seems difficult. Our numerical tests are performed on the same equations considered in \cite{han2018solving} which have benchmarked solutions derived at single points by trusted modern numerical methods. Other aspects that we investigate are numerical tests of the average approximation error for the solution in {\it whole domains}, time scaling and complexity of the algorithm. Further systematic explorations of even higher order discretizations, as well as their implementation on further algorithmic improvements proposed in the literature, is beyond the scope of this paper.

The paper is organized as follows. In Sec. 2 we introduce standard background material: backward Kolmogorov equations, the associated FBSDE's and the non-linear Feynman-Kac formula. In Sec. 3 we apply the formalism with two different 
discretizations, namely the Milstein and Leimkhuler and Matthews schemes. 
In Sec. 4 we describe for each discretization the associated architecture of the deep network and learning algorithm. Finally, in Sec. 5, we discuss numerical results for high dimensional Black-Scholes, an Allen-Cahn, non-linear exactly solvable diffusion, and Hamilton-Jacobi-Bellman equations. All the codes can be downloaded at \cite{Marino2019GITHUB}.

\section{Non-linear Kolmogorov equations and variational formulation of their solution}\label{backward-K-equation}

We give standard background material on Kolmogorov's equation, the associated forward-backward SDEs, and non-linear Feynman-Kac formula. This later formula can be seen as the solution of a variational problem which is at the basis of the deep learning algorithm. Following \cite{han2018solving} we give the formulation in terms of backward Kolmogorov equations where the {\it terminal} condition is given. A large class of semilinear parabolic forward PDEs with given {\it initial} condition can be transformed into backward Kolmogorov equations considered below (by a simple time reversal transformation). 


\subsection{Backward Kolmogorov equation and stochastic differential equations}
\label{subsec::Keq}

Let $g : \mathbb{R}^d \times [0, T] \to \mathbb{R}$ a function which is of class $C^2$ in $\vec x$ and $C^1$ in $t$. 
Let $A:\mathbb{R}^d \times [0, T] \to  \mathbb{R}^d$ and $B: \mathbb{R}^d\times [0,T] \to  \mathbb{R}^{d \times d}$ be piecewise differentiable (in $\vec x)$  vector and tensor transport coefficients. Let also 
$f:[0,T]\times \mathbb{R}^d\times \mathbb{R}_+\times \mathbb{R}^d \to \mathbb{R}$ be a known function which captures non-linearities.
Non-linear backward Kolmogorov equations \cite{kolmogoroff1933theorie} describe a time-evolution of the form
\begin{align}
\label{Kolmogoroveq}
\frac{\partial g(\vec{x}, t)}{\partial t}= & -A(\vec{x}, t)\cdot \nabla_{\vec{x}} g(\vec{x}, t) - \frac{1}{2}
\text{Tr} \big(BB^{T}(\vec{x},t)\text{Hess}_{\vec{x}} g(\vec{x}, t)\big) 
\nonumber \\ &
- f\big(t, \vec{x}, g(\vec{x}, t), B^{T}(\vec{x},t)\nabla g(\vec{x}, t)\big), \qquad t<T
\end{align} 
given some specific {\it terminal} condition $g(\vec{x}, T)=\phi(\vec{x})$ where $\phi: \mathbb{R}^d \to \mathbb{R}$ is a {\it known} function. Here $B^{T}$ denotes the transpose of matrix $B$, $\nabla_{\vec x} g(\vec{x}, t)$ and $\text{Hess}_{\vec{x}} g(\vec{x}, t) = 
\big(\frac{\partial^2 g}{\partial x_i\partial x_j}\big)_{i,j=1}^d$ denote the gradient and the Hessian of $g$ with respect to $\vec{x}$ respectively, \text{Tr} denotes the trace of a matrix. To fix ideas we consider
real numbers $x_1, x_2 \in \mathbb{R}$ with $x_1<x_2$ and we suppose that our goal is to compute the function on the whole cube $ \vec{x} \in  [x_1, x_2]^d  \mapsto \mathbb{R}$ for $t<T$.


It is well known that when the non-linearity is absent, i.e., $f=0$, the adjoint of \eqref{Kolmogoroveq} describes the time evolution of a  probability density of an underlying $d$-dimensional stochastic process $\{\vec{X}(t), t\in [0,T]\}$ given by the (forward) SDE
\cite{kolmogoroff1933theorie},
\begin{equation}
    \label{eq:Langevin}
    d\vec{X}(t) = A[\vec{X}(t), t]dt + B[\vec{X}(t), t]d\vec{W}(t).
\end{equation}
where $\{\vec{W}(t), t\in [0,T]\}$ is the $d$-dimensional Wiener process, and we have $g[\vec X(T), T] = \phi[\vec X(T)]$.

Consider now stochastic trajectories generated by \eqref{eq:Langevin} (with say the initial condition $\vec{X}(0) = \vec x$). The solution of the
non-linear Kolmogorov equation evaluated along a stochastic trajectory, i.e., $g(\vec{X}(t), t)$ satisfies a (backward) SDE which is at the basis of the deep learning algorithms discussed later on. Applying Ito's chain rule on $g[\vec{X}(t), t]$,
\begin{align}
    dg[\vec{X}(t),t]  =  \frac{\partial g}{\partial t}[\vec{X}(t), t] dt + \nabla_{\vec{x}} g[\vec{X}(t),t]^T d\vec X(t)+ \frac{1}{2} \text{Tr} \big(BB^{T}[\vec{X}(t),t]\text{Hess}_{\vec{x}} g(\vec{X}(t), t)\big)dt
\end{align}
and then using equation \eqref{Kolmogoroveq} we obtain
\begin{align}
    \label{d_sol}
    dg[\vec{X}(t),t] =
    - f\big(t, \vec{X}(t), g[\vec{X}(t),t], B[\vec{X}(t),t]^T\nabla_{\vec{x}}g[\vec{X}(t),t]\big) dt
    +
    \nabla_{\vec{x}} g[\vec{X}(t),t]^T B[\vec{X}(t),t] d\vec{W}(t) .
\end{align}
Integrating this equation on an interval $[t_1,t_2]\subset [0, T]$ we get
\begin{align}
    \label{sol}
    g[\vec{X}(t_2),t_2] = & \, g[\vec{X}(t_1),t_1] -\int_{t_1}^{t_2}  f\big(t, \vec{X}(t), g[\vec{X}(t),t], B[\vec{X}(t),t]^T \nabla_{\vec{x}} g[\vec{X}(t),t]\big)dt
    \nonumber \\ & 
 + \int_{t_1}^{t_2}(\nabla_{\vec{x}}g[\vec{X}(t),t])^T B[\vec{X}(t),t]d\vec{W}(t)     . 
\end{align}
This is the backward SDE. As will be seen in Sect. 4, the various discretized forms of equations \eqref{eq:Langevin} and \eqref{sol} constitute the basis of the deep networks. 

When \eqref{sol} is applied to $t_1=0$ and $t_2 =T$ for paths conditioned to satisfy $\vec{X}(0) = \vec{x}$ and the terminal condition $g(\vec{X}(T), T) = \phi(\vec{X}(T))$ we get the so-called non-linear Feynman-Kac formula
\begin{align}
    \label{intermediate}
    g(\vec{x},0) = & \, \phi(\vec{X}(T)) +\int_{0}^{T}  f\big(t, \vec{X}(t), g[\vec{X}(t),t], B[\vec{X}(t),t]^T \nabla_{\vec{x}} g[\vec{X}(t),t]\big)dt
    \nonumber \\ & 
 - \int_{0}^{T}(\nabla_{\vec{x}}g[\vec{X}(t),t])^T B[\vec{X}(t),t]d\vec{W}(t)     . 
\end{align}
Taking the expectation $\mathbb{E}_{\vec x}$ over the stochastic process $\{\vec{X}(t), 0\leq t \leq T\}$ conditioned on $\vec X(t) = \vec x$ we obtain the following form of the non-linear Feynman-Kac formula
\begin{align}\label{eq:non-lin-FK}
    g(\vec x, 0) = \mathbb{E}_{\vec x}\phi[\vec X(T)] + \int_{0}^{T}  \mathbb{E}_{\vec x}f\big(s, \vec{X}(s), g[\vec{X}(s),s], B[\vec{X}(s),s]^T \nabla_{\vec{x}} g[\vec{X}(s),s]\big)ds.
\end{align}
Note that the expectation coming from  the second term on the r.h.s of \eqref{intermediate} vanishes because \\
$
(\nabla_{\vec{x}}g[\vec{X}(s),s])^T B[\vec{X}(s),s]
$
is independent of the increment $d\vec W(s)$ (i.e., it is a non-anticipating functional 
\cite{gardiner1985handbook}).
In the linear case $f=0$ equation \eqref{eq:non-lin-FK} is nothing else than the usual Feynman-Kac formula. 

\subsection{Variational formulation}\label{sec:varform}

Let us explain the main idea behind the deep learning algorithms introduced in \cite{han2018solving} which is based on a variational reformulation of the non-linear Feynman-Kac formula \eqref{intermediate}.
Consider a variational class of functions or "deep networks" $\mathcal{NN}_{g_0}(\vec x\mid \vec\theta_0) \in \mathbb{R}$ and $\mathcal{NN}_{\nabla g_s}(\vec x\mid \vec\theta(s)) \in \mathbb{R}^d$ where $\vec\theta_0$ and $\{\vec\theta(s)$, $0\leq s\leq T$\} is a set of "weights" to be optimized. We want to minimize the  loss functional
\begin{align}
\mathcal{L}[\vec\theta_0, \vec\theta(\cdot)] = & \mathbb{E}\Big[ \Big\vert \phi(\vec X(T)) -   \Big\{ \mathcal{NN}_{g_0}(\vec x\vert \theta_0)   
-\int_0^T ds f(s, \vec X(s), g(\vec X(s), s),  \mathcal{NN}_{\nabla g_s}(\vec X(s) \vert \vec\theta(s)))   
\nonumber \\ &
+ \int_0^T (\mathcal{NN}_{\nabla g_s}(\vec x\vert \vec\theta(s)))^T  d\vec W(s) \Big\}\Big\vert^2     \Big]
\label{loss-cont-language}
\end{align}
where $\vec X(s)$ and $g(\vec X(s),s)$ are solutions of the SDE's
\begin{align}\label{sde-var}
\begin{cases}
    d\vec X(s) = A(X(s),s) ds + B(X(s), s)^Td \vec W(s),
    \\
    dg[\vec{X}(s),s] \!=\!
    - f\big(s, \vec{X}(s), g[\vec{X}(s),s], \mathcal{NN}_{\nabla g_s}(\vec X(s)\vert \vec\theta(s))\big) ds
    \!+\!
    (\mathcal{NN}_{\nabla g_s}(\vec X(s)\vert \vec\theta(s)))^T \!d\vec{W}(s) 
    \end{cases}
\end{align}
with $g(\vec X(0), 0) = \mathcal{NN}_{g_0}(x\vert \vec\theta_0)$. The expectation in \eqref{loss-cont-language} is over the Brownian trajectories $\vec W(s)$, $0\leq s\leq T$ and over $\vec X(0)\in D$ uniformly distributed in a domain $D \subset \mathbb{R}^d$.  

The non-linear Feynman-Kac formula \eqref{intermediate} tells us that as long as the class of the deep networks is sufficiently expressive, 
there exists a minimizer $\widehat\theta_0$, $\widehat\theta(\cdot)$ such that 
once the loss is sufficiently small, the approximation of  the deep networks to the PDEs solutions becomes more accurate
, and $\mathcal{NN}_{g_0}(\vec x\vert \widehat\theta_0) \approx g(\vec x\vert 0)$ and $\mathcal{NN}_{\nabla g_s}(\vec x\vert \widehat\theta(s)) \approx B(\vec x, s)^T\nabla_{\vec x} g(\vec x, s)$.\footnote{
A possible route to make this argument rigorous would be to use methods of \cite{han2020convergence, jiang2021convergence}. Rigorous aspects are beyond the scope of this paper and we do not further address this point here.}
The neural networks investigated in section \ref{sec:neural} are based on discretization schemes of \eqref{loss-cont-language} and \eqref{sde-var}.

Finally we wish to remark here that the minimizers satisfy a formula analogous to \eqref{eq:non-lin-FK}. 
The argument is rather standard but we give the details for the convenience of the reader.
The loss functional is of the form 
\begin{align}\label{def0}
\mathcal{L}[\theta_0, \theta(\cdot)] = \mathbb{E}[\{\mathcal{NN}_{g_0}(\vec x\mid \theta_0) - \mathcal{F}(\vec W(\cdot) \mid \theta_0, \theta(\cdot))\}^2]
\end{align}
where $\mathcal{F}$ is a complicated functional of the Brownian paths $\vec W(s)$, $0\leq s\leq T$. We have
\begin{align}
\mathcal{L}[\theta_0, \theta(\cdot)] = & \mathbb{E}[\{\mathcal{NN}_{g_0}(\vec x\mid \theta_0) - \mathbb{E}_{\vec x}[\mathcal{F}(\vec W(\cdot) \mid \theta_0, \theta(\cdot))]\}^2]
\nonumber \\ &
+ \mathbb{E}[\{\mathbb{E}_{\vec x}[\mathcal{F}(\vec W(\cdot) \mid \theta_0, \theta(\cdot))] - \mathcal{F}(\vec W(\cdot) \mid \theta_0, \theta(\cdot))\}^2]
\nonumber \\ &
+ 2 \mathbb{E}[ \{\mathbb{E}_{\vec x}[\mathcal{F}(\vec W(\cdot) \mid \theta_0, \theta(\cdot))] - \mathcal{F}(\vec W(\cdot) \mid \theta_0, \theta(\cdot))\}
\{\mathcal{NN}_{g_0}(\vec x\mid \theta_0) - \mathbb{E}_{\vec x}[\mathcal{F}(\vec W(\cdot)\mid \theta_0, \theta(\cdot))]\} ]
\end{align}
where we recall that $\mathbb{E}_{\vec x}$ is the expectation conditioned on $\vec X(0) = \vec x$ (in other words the expectation only over $\vec W(s)$). Since $\mathbb{E} = \mathbb{E}_{\vec X(0)} \mathbb{E}_{\vec x}$ we find that 
the first bracket $\{\cdots \}$ in third term vanishes. Therefore
\begin{align}\label{an-inequality}
\mathcal{L}[\theta_0, \theta(\cdot)] \geq 
\mathbb{E}[\{\mathbb{E}_{\vec x}[\mathcal{F}(\vec W(\cdot) \mid \theta_0, \theta(\cdot))] - \mathcal{F}(\vec W(\cdot) \mid \theta_0, \theta(\cdot))\}^2]
\end{align}
Comparing \eqref{def0} with inequality \eqref{an-inequality} we see that this inequality is saturated for a choice of the neural network parameters 
$\widehat\theta_0$, $\widehat\theta(\cdot)$ such that
\begin{align}
\mathcal{NN}_{g_0}(\vec x \vert \widehat\theta_0) 
& \approx \mathbb{E}[\mathcal{F}(\vec W(\cdot)\vert \widehat\theta_0, \widehat\theta(\cdot))]
\nonumber \\ &
= 
\mathbb{E}[\phi(X(T)] + \int_0^T ds \mathbb{E}[f(s, \vec X(s), g(\vec X(s), s),  \mathcal{NN}_{\nabla g_s}(\vec X(s) \vert \widehat\theta(s)))]
\label{approx-FK}
\end{align}
We thus arrive at the conclusion that the loss is minimized when the neural network satisfies \eqref{approx-FK}.
In the linear case $f=0$ we recover again the usual Feynman-Kac formula.

\section{Stochastic Numerical Integration}
\label{sec::SNI}
In this section we present three discretizations of the SDE \eqref{eq:Langevin}. The first one, the Euler-Maruyama scheme, is the simplest and the one considered in \cite{han2018solving}. It is in fact the first of a hierarchy of higher-order discretizations with smaller errors \cite{kloeden2013numerical}, which are however often difficult to apply in multidimensional settings. Here we discuss the next best discretization in this hierarchy, namely the second-order Milstein scheme.  The third discretization that we consider is a first order scheme recently proposed by Leimhkhuler and Matthews in the context of high-friction limits of Langevin dynamics used for sampling \cite{leimkuhler2012rational}. 

Let us first explain the general principle of the Euler-Maruyama and Milstein schemes. 
Integrating \eqref{eq:Langevin} over a small time interval $[t_1, t_2]$ and applying Ito's chain rule lemma on
$A[\vec{X}(t), t]$ and $B[\vec{X}(t), t]$ one obtains (with component-wise notation)
\begin{align}
\label{eq_X_ito_transformed}
\begin{split}
&X_i(t_2)=X_i(t_1)+A_i[\vec{X}(t_1),t_1] (t_2-t_1) + \sum_{j=1}^d B_{ij}[\vec{X}(t_1),t_1]  (W_j(t_2) - W(t_1))+\\
&+\sum_{j,k,l=1}^d B_{kl}[\vec{X}(t_1),t_1] \partial_{x_k} B_{ij}[\vec{X}(t_1),t_1]\int_{t_1}^{t_2}dW_j(s)\int_{t_1}^{s}  dW_l(s')
+\cdots
\end{split}
\end{align}
In developing equation \eqref{eq_X_ito_transformed} we have kept all terms "linear in $\delta t = t_2-t_1$" and have discarded all higher terms. 

For implementing a discretization of the BSDE we divide the time interval $[0,T]$ into $N$ sub-intervals of size $\tau=T/N$ at points $\tau_n=n\tau$, so that a function can be evaluated at times 
\begin{equation}
    \label{timediscrete}
    \tau_0=0, \, \tau_1, \, \tau_2, \, \tau_3, \dots \, \tau_{N-1}, \, \tau_N=T.     
\end{equation}
Thus, the corresponding Wiener increments become $\Delta \vec{W}^n=\vec{W}(\tau_{n+1})-\vec{W}(\tau_{n})$.  The (exact) stochastic process at the same times is denoted as $\vec{X}^n=\vec{X}(\tau_n)$, and satisfies (from equation \eqref{eq_X_ito_transformed})
\begin{align}
\label{exact_eq_transf}
X^{n+1}_i= & X^n_i+A_i[\vec{X}^n,\tau_n]\tau + \sum_{j=1}^d B_{ij}[\vec{X}^n,\tau_n]\Delta W^n_j 
\nonumber \\ &
+ \sum_{j,k,l=1}^d B_{kl}[\vec{X}^n,\tau_n] \partial_{x_k} B_{ij}[\vec{X}^n,\tau_n]\mathcal{I}_{jl}^n +\cdots
\end{align}
where 
$$
\mathcal{I}_{jl}^n=\int_{\tau_n}^{\tau_{n+1}}dW_j(s)\int_{\tau_n}^{s}  dW_l(s').
$$
This multiple stochastic integral cannot be simply expressed in terms of the increments $\Delta {W}^n_{l}$ and $\Delta {W}^n_{j}$ of the components of the Wiener process. However, an exercise in Ito calculus shows that the symmetrised quantity has a simple expression
$$
\mathcal{I}_{jl}^n + \mathcal{I}_{lj}^n = \Delta W_j^n \Delta W_l^n - \tau\delta_{jl},
$$
and therefore, when the tensor $\mathbf{B}$ satisfies the so-called "commutativity" condition \cite{glasserman2004monte}:
\begin{equation}
\label{commutative_noise}
\sum_{k=1}^d B_{kj}[\vec{X}^n,\tau] (\partial_{x_k} B_{il}[\vec{X}^n,\tau])=\sum_{k=1}^d B_{kl}[\vec{X}^n,\tau] (\partial_{x_k} B_{ij}[\vec{X}^n,\tau]),
\end{equation}
we can replace $\mathcal{I}_{lj}^n$ in \eqref{exact_eq_transf} by $\frac{1}{2}(\Delta W_j^n \Delta W_l^n - \tau\delta_{jl})$. This leads to
\begin{align}
    X^{n+1}_i= & X^n_i+A_i[\vec{X}^n,\tau_n]\tau + \sum_{j=1}^d B_{ij}[\vec{X}^n,\tau_n]\Delta W^n_j 
\nonumber \\ &
+ \frac{1}{2}\sum_{j,k,l=1}^d B_{kl}[\vec{X}^n,\tau_n] (\partial_{x_k} B_{ij}[\vec{X}^n,\tau_n])(\Delta W_j^n \Delta W_l^n - \tau\delta_{jl}) +\cdots
\label{ito-taylor}
\end{align}
This Taylor-Ito expansion forms the basis of the Euler-Maruyama and Milstein discretizations discussed below. 

It is intriguing  to understand when the commutativity condition holds. Equation \eqref{commutative_noise}, following the words of Gardiner in \cite{gardiner1985handbook}, \textit{"is not uncommonly satisfied"}. Gardiner in \cite{gardiner1985handbook} lists the notable cases where the commutative condition holds. However,  his definition of the matrix $B$ is much general then the one chosen by us, and his list is bigger than our. We confine ourselves in the case where $B: \mathbb{R}^d\times [0,T] \to  \mathbb{R}^{d \times d}$.	Under this assumption the most interesting cases, where the commutativity conditions holds, are: i) the case when we have additive noise, i.e. when the coefficients $B_{ij}$  may depend explicitly by time only. ii) The case when $B_{ij}(\vec{x}(t), t)=x_iG_{i,j}$, where $G$  is an arbitrary constant square $d \times d$ matrix. In this particular case, eq. \eqref{commutative_noise} is satisfied. Indeed, 
$\sum_{k} B_{kl} (\partial_{x_k} B_{ij})= \sum_{k} x_kG_{kl}\delta_{ik}G_{ij}= x_iG_{ij}G_{il}=\sum_{k} B_{kj}(\partial_{x_k} B_{il})$. This calculation can be generalized for any elementary scalar continuous function  $f: \mathbb{R} \to \mathbb{R}$ such that   $B_{ij}(\vec{x}(t), t)=f(x_i)G_{i,j}$  where $G$  is an arbitrary constant square $d \times d$ matrix. Indeed, $\sum_{k} B_{kl} (\partial_{x_k} B_{ij})= \sum_{k} f(x_k)G_{kl}f'(x_i)\delta_{ik}G_{ij}= f(x_i)G_{ij}f'(x_i)G_{il}=\sum_{k} B_{kj}(\partial_{x_k} B_{il})$, which proofs that eq. \eqref{commutative_noise} is satisfied.
iii) It is a special case of ii), i.e., when $G$ is a $d \times d$ matrix and $G_{ij} = c \delta_{ij}$, where $c$ is a constant \cite{gardiner1985handbook}. This case has been studied in detail in \cite{kloeden2013numerical}.

The most notable case for us is iii. We observe, indeed, that in finance, physics and in many other disciplines this case is the most promising. Indeed, in many real examples, the diffusion tensors $D=BB^T$ is diagonal or it is chosen to be real symmetric, or Hermitian, or normal, and, thus, diagonalizable by orthogonal, for real symmetric, or unitary, for the other two cases, matrices. 

	For example, in finance there are securities which are traded in the stock market, and have a value which fluctuates with   time depending on market conditions. If the value of one item of stock as traded on the stock market is $x(t)$, the appropriate equation for the time dependence on this value is written as a stochastic differential equation with a coefficient in front of the Wiener process equal to $\sigma (t)x(t)$, where $\sigma (t)$ is the volatility.  $d$ independent securities, then, will have the $B$ coefficient diagonal, as in case iii, satisfying the commutative condition. 
	
	In physics, instead, the  diffusion of micron-sized symmetrical particles, as spheres, ellipsoids, and bodies of revolution with three mutually perpendicular symmetry planes \cite{brenner1967coupling}, settling in a gravitational field is an intriguing phenomenon. 
	This diffusion process is studied for understanding sedimentation and aggregation during sedimentation processes \cite{brenner1981taylor, brenner1981taylor2, marino2016advective}. In a sedimentation process, particles can traslate and rotate during their motion, and each particle can be described by $d=6$ degrees of freedom, $3$ for rotation and $3$ for translation. In a laboratory frame, the tensors, which are defined positive matrices describing the friction for the rotational motion and the translational one, depend on the parameters specifying the particle orientation. They couple, therefore, the two motions. However, for a symmetrical particle there exists always a unique geometrical point, which is called \textit{center of diffusion}, where the two friction tensors are diagonalizble \cite{brenner1967coupling}.  This center of diffusion point allows to make a transformation of coordinates where the two friction tensors are diagonal in the new basis. Collecting the two diagonal friction tensors in a new $d=6$ dimensional space, one can build examples where the temperature in the new coordinates allows to have the $6 \times 6$  diagonal matrix $B$ of the form given by iii. For more information about rotational and translational diffusion processes, we refer the interesting reader to \cite{brenner1967coupling, brenner1981taylor, brenner1981taylor2, marino2016advective, marino2016entropy, aurell2016diffusion} and references therein.


\subsection{The Euler-Maruyama scheme}

This is the simplest and most used approximation for computing the solution of a stochastic differential equation. The discrete stochastic process $\vec{Y}^n_{EM}$ is defined by
\begin{align}
\label{discretiz_EM}
(Y_{EM})^{n+1}_i=(Y_{EM})^n_i+A_i[\vec{Y}^n_{EM},\tau_n]\tau + \sum_{j=1}^d B_{ij}[\vec{Y}^n_{EM},\tau_n]\Delta W^n_j 
\end{align}
This scheme obviously corresponds to the truncated Ito-Taylor
expansion retaining only terms of order one in \eqref{ito-taylor} (see \cite{kloeden2013numerical}). 
This discretization has a "strong error", defined as the root-mean-square of the difference with the exact solution, of order $O(\sqrt \tau)$ after a finite time $n = O(N)$ \cite{gardiner1985handbook, kloeden2013numerical}.

\subsection{The Milstein scheme}
\label{Milsteinalgo}
The Milstein scheme represents a high-order approximation of the solution of the BSDE, which keeps in \eqref{exact_eq_transf} all multiple stochastic integrals proportional to $O(\tau)$ (we recall that $\Delta \vec{W}^n \sim \sqrt{\tau}$). 
In this paper we only investigate cases where the relation \eqref{commutative_noise} is valid and therefore we can consider the expansion 
\eqref{ito-taylor}.
Denoting by $\vec{Y}^n_{M}$ the discrete stochastic process of the Milstein scheme \cite{kloeden2013numerical, gardiner1985handbook, jentzen2015milstein}, we have:
\begin{equation}
\label{discretiz_M}
\begin{split}
&(Y_{M})^{n+1}_i=(Y_{M})^n_i+A_i[\vec{Y}^n_M,\tau_n]\tau + \sum_{j=1}^d B_{ij}[\vec{Y}^n_M,\tau_n]\Delta W^n_j +\\ &+\frac{1}{2}\sum_{j,k,l=1}^d B_{kl}[\vec{Y}^n_M,\tau_n] (\partial_{x_{k}}B_{ij}[\vec{Y}^n_M,\tau_n]) (\Delta W^{n}_{j}\Delta W^{n}_{l}-\delta_{jl}\tau)
\end{split}
\end{equation}
This scheme has a strong error $O(\tau)$ after a finite time $n=O(N)$ \cite{gardiner1985handbook, kloeden2013numerical}.

\subsection{The Leimkuhler-Matthews scheme}
\label{sec::LM}
The Leimkuhler-Matthews discretization is a simple modification of the Euler-Maruyama approximation involving a non-Markovian (colored noise) random process. This discretization is derived in \cite{ leimkuhler2012rational} by looking at various approximations for the high-friction limit of the Langevin dynamics (overdamped dynamics), and it is shown that the statistical error improves to $O(\tau^2)$. Although this method has been developed for molecular sampling, we choose it among others because it is linear in the Brownian increments, and it has the same computational complexity of the two methods above described, but involving a non-Markovian (colored noise) random process.

The Leimkuhler-Matthews discrete stochastic process is
\begin{align}
\label{discretiz_LM}
(Y_{LM})^{n+1}_i=(Y_{LM})^n_i+A_i[\vec{Y}^n_{LM},\tau_n]\tau + \frac{1}{2}\sum_{j=1}^d B_{ij}[\vec{Y}^n_{LM},\tau_n](\Delta W_j^n +\Delta W_j^{n+1}). 
\end{align}
This method samples  better trajectories of overdamped Langevin dynamics when $\tau$ is small enough. However we will observe that the deep network based on it does not return a better results as compared to the Euler-Maruyama and Milstein discretizations. This is the case even when the diffusion tensor is independent of space (in which case the Euler-Maruyama and Milstein schemes are equivalent) and seems to be related to the difficulty of taking a small enough $\tau$ when the deep networks are implemented.

\section{Neural Networks for solving backward Kolmogorov's PDEs}
\label{sec:neural}

The goal is to determine an approximation  $g_0(\vec x\vert \vec{\theta}_{g_0})$ for the backward solution $g(\vec x,0)$ of the non-linear Kolmogorov equation \eqref{Kolmogoroveq} when the terminal condition is known $g(\vec x, T)=\phi(\vec x)$. Here $\vec{\theta}_{g_0}$ are parameters that have to be learned. We use the main idea put forward in \cite{han2018solving}.
The function $g_0(\vec x\vert \vec{\theta}_0)$ is determined by a neural network which implements a discrete version of equation \eqref{sol}. The network takes as {\it input} a discrete stochastic trajectory 
$$
(\vec{Y}^0, \Delta \vec{W}^0), (\vec{Y}^1,\Delta\vec{W}^1), \cdots ,(\vec{Y}^N, \Delta \vec{W}^N)
$$ 
and computes approximations for $g(\vec{Y}^0,0), g(\vec{Y}^1, \tau_1), \cdots, g(\vec{Y}^N, \tau_N)$. The last approximation for 
$g(\vec{Y}^N, \tau_N)$ is the final {\it output} of the network
\begin{align}
\label{estimate}
\hat{g}(\{\vec{Y}^{n}, \Delta\vec{W}^n\}_{0\leq n \leq N} \mid \vec{\Theta}) \approx g(\vec{Y}^N, \tau_N)
\end{align}
and depends on the whole discrete trajectory as well as a collection $\vec{\Theta}$ of all parameters that have to be learned. Here learning amounts to minimize the following quadratic loss
\begin{equation}
\label{loss_function}
\mathcal{L}(\mathbf{\Theta}):=\mathbb{E}\left[|\phi(\vec{Y}^{N})-
\hat{g}(\{\vec{Y}^{n}\}_{0\leq n \leq N}, \{\Delta \vec{W}^{n}\}_{0\leq n \leq N}\vert \vec{\Theta})|^2 \right].
\end{equation}
Each discretization of the BSDE leads to a different discretization for equation \eqref{sol} which in turn leads to a different architecture for the network. In the next paragraphs we present the details of the three architectures corresponding the Euler-Maruyama, Milstein and Leimkhuler-Matthews discretizations.  We first repeat the first architecture found in \cite{han2018solving} for pedagogic reasons, while the other two are slightly more complicated variants.

\begin{figure}[h!]
    \centering
    \includegraphics[width=0.9\columnwidth]{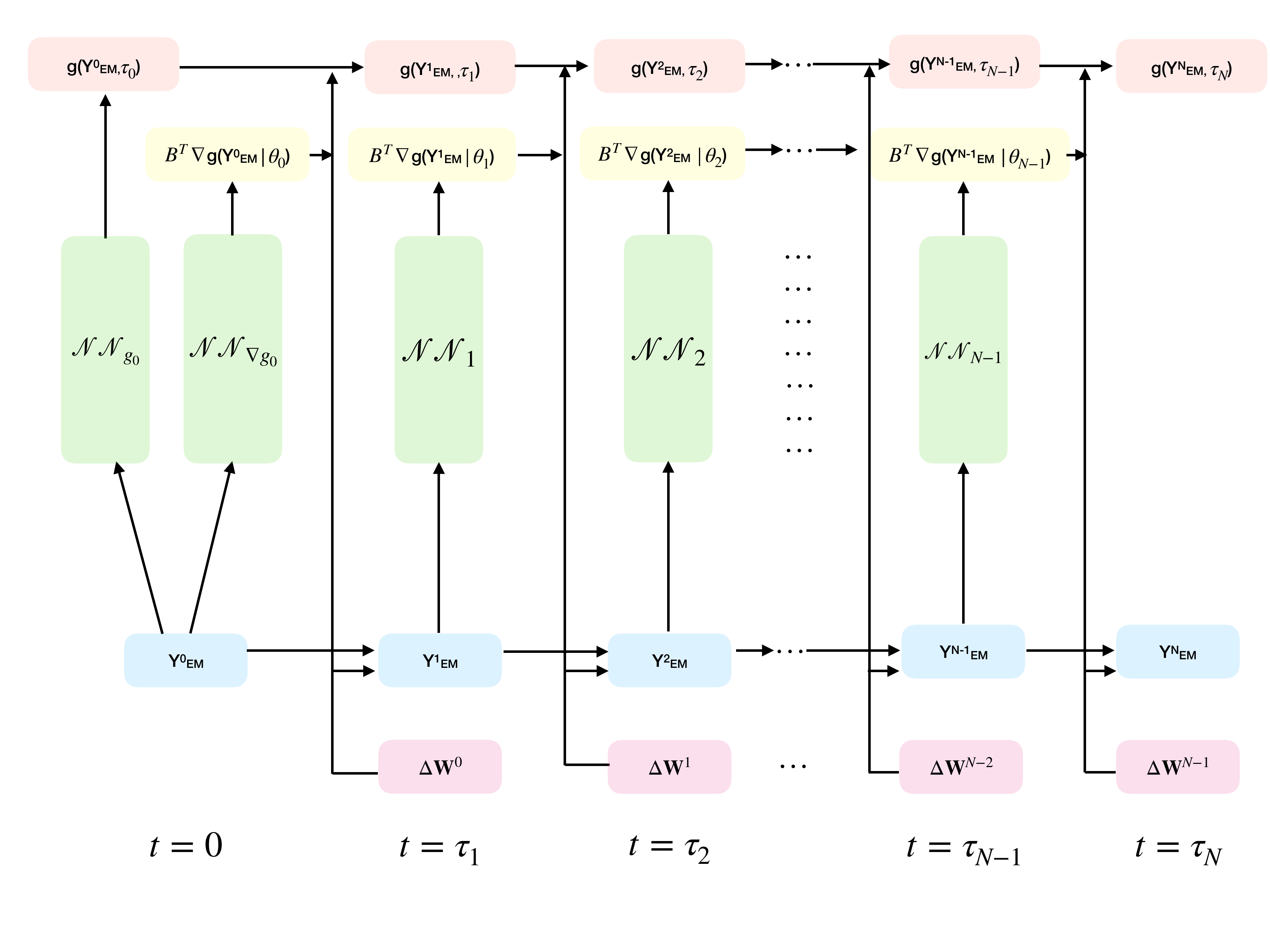}
    \caption{\small $\mathcal{DNN}_{EM}$: the network obtained from the Euler-Maruyama scheme. Time flows horizontally, from left to right, while vertically, the sub-networks output approximations of $g$ and $\nabla g^T B$ (and its transpose) at each time step. $\mathcal{DNN}_{EM}$ takes as {\it input} discrete stochastic trajectories $\{\vec{Y}^{n}_{EM}, \Delta{\vec W}^n\}_{0\leq n \leq N}$ and {\it outputs} the value of $\hat{g}(\{\vec{Y}^{n}_{EM}, \Delta{\vec W}^n\}_{0\leq n \leq N}\vert \vec{\Theta})$. Specifically, at time $t=0$, the neural network takes as input the set of points of the stochastic trajectories at time $t=0$. The data are processed by two separate sub-neural networks $\mathcal{NN}_{g_{0}}$ and $\mathcal{NN}_{\nabla g_{0}}$ . Their output allows to compute the approximate value of $g(\vec{Y}^{1}_{EM},\tau_1)$ using the recursive equation \eqref{sol_BSDE_discret}.
     At this point time flows and at $t=\tau_1$ and $(\vec{Y}^{1}_{EM}, \Delta{\vec W}^1)$ are given to the sub-network $\mathcal{NN}_1$. This sub-network returns $(\nabla g^T B)[\vec{Y}^1_{EM}, \tau^1]$ (see \eqref{_grad_approx}) 
     used to get the approximate value of $g(\vec{Y}^{2}_{EM},\tau_2)$. This procedure continues till time $\tau_N$ where the approximate value for $g(\vec{Y}^{N}_{EM},\tau_N)$ is obtained (see \eqref{estimate}).}
    \label{fig:DNN_EM}
\end{figure}

\subsection{Architecture based on Euler-Maruyama scheme}\label{subsec::EM}
Using the Euler-Maruyama scheme \eqref{discretiz_EM} the discrete form of equation \eqref{sol} becomes
\begin{align}
    \label{sol_BSDE_discret}
    & g[\vec{Y}^{n+1}_{EM},\tau_{n+1}] = 
    g[\vec{Y}^n_{EM},\tau_n] + \nabla_{\vec x} g[\vec{Y}^n_{EM},\tau_n]^T B[\vec{Y}^n_{EM},\tau_n]\Delta \vec{W}^n
    \nonumber \\ &
    -f(\tau_n, \vec{Y}^n_{EM}, g[\vec{Y}^n_{EM},\tau_n], B[\vec{Y}^n_{EM},\tau_n]^T \nabla_{\vec x} g[\vec{Y}^n_{EM},\tau_n]) \tau 
\end{align}
A key observation here is that once $g[\vec{Y}^n_{EM},\tau_n]$ is determined, the only term which is unknown in \eqref{sol_BSDE_discret} is  $\nabla g[{Y}^n_{EM},\tau_n]^T B[\vec{Y}^n_{EM},\tau_n]$ (and its transpose). Such a gradient is a function of the path $\vec{Y}^n_{EM}$ at time $\tau_n$ and can be approximated by a feed-forward neural network  $\mathcal{N}\mathcal{N}_n$ whose output we denote
\begin{equation}
\label{_grad_approx}
 (\nabla g^T B)[{\vec Y}^n_{EM}\vert \vec{\theta}_n] \approx \nabla_{\vec x} g[{\vec Y}^n_{EM},\tau_n]^T B[\vec{Y}^n_{EM},\tau_n]
\end{equation}
for $n=1,\dots,N-1$, where $\vec{\theta}_n$ are parameters to be learned. We assemble all these networks 
$\mathcal{N}\mathcal{N}_1, \cdots, \mathcal{N}\mathcal{N}_{N-1}$
into one big deep neural network, called $\mathcal{DNN}_{EM}$, adding also two extra feed-forward networks $\mathcal{NN}_{g_{0}}$ and $\mathcal{NN}_{\nabla g_{0}}$ at time $t=0$ needed to approximate the initial values $g[\vec{Y}^0_{EM}, 0]$ and 
$\nabla g[{\vec Y}^0_{EM},0]^T B[\vec{Y}^0_{EM},0]$. The complete architecture as well as the sequence of operations are explained in figure \ref{fig:DNN_EM}. $\mathcal{DNN}_{EM}$ takes as {\it input} discrete stochastic trajectories $\{\vec{Y}^{n}_{EM}, \Delta{\vec W}^n\}_{0\leq n \leq N}$
and {\it outputs} the value of $\hat{g}(\{\vec{Y}^{n}_{EM}, \Delta{\vec W}^n\}_{0\leq n \leq N}\vert \vec{\Theta})$
where $\vec{\Theta}=\{\vec{\theta}_{g_{0}}, \vec{\theta}_{\nabla g_{0}}, \vec{\theta}_1, \dots, \vec{\theta}_{N-1}\}$ is the total set of parameters to be learned. As indicated above learning amounts to minimize the loss function \eqref{loss_function} which can be done by standard stochastic gradient descent (SGD) algorithms. In this process the initial network $\mathcal{N}\mathcal{N}_{g_0}$ is also learned and allows to compute $g_0(\vec{x}\vert\theta_0) \approx g(\vec{x}, 0)$.

In this work we choose the same internal architecture as in \cite{han2018solving} which is simple and has proven effective.
$\mathcal{NN}_{g_{0}}$ has 4 layers. The first input layer has $d$ neurons to input $\vec{Y}^0_{EM}$, the second and third hidden layers both have $d+10$ neurons, and the fourth output layer has one neuron to return the approximate value of $g(\vec{Y}^0, 0)$. The activation functions are simply given by ${\rm ReLu}$'s.
The sub-networks $\mathcal{NN}_{\nabla g_{0}}$ and $\mathcal{N}\mathcal{N}_n$ have the same structure except for the output layer that must contain $d$ neurons to output \eqref{_grad_approx}.

\subsection{Architecture based on Milstein scheme}\label{Milstein_DNN}
Using the higher-order Milstein discretization \eqref{discretiz_M} the discrete form of equation \eqref{sol} becomes
\begin{align}
    \label{sol_BSDE_discret_MILSTEIN}
    g[\vec{Y}^{n+1}_{M}, & \tau_{n+1}] =  g[\vec{Y}^n_{M},\tau_n]
    \nonumber \\ &
    -f(\tau_n, \vec{Y}^n_{M}, g[\vec{Y}^n_{M},\tau_n], B[\vec{Y}^n_{M},\tau_n]^T \nabla_{\vec x} g[\vec{Y}^n_{M},\tau_n])\tau
    \nonumber \\ & 
    + \nabla_{\vec x}g[\vec{Y}^n_{M},\tau_n]^T B[\vec{Y}^n_{M},\tau_n]\Delta \vec{W}^n
    \nonumber \\ & 
    +\frac{1}{2}\sum_{i, j=1}^d \Big\{\sum_{l,k=1}^d(\partial_{x_l}g[\vec{Y}^n_{M},\tau_n]) B_{kl}[\vec{Y}^n_{M},\tau_n] (\partial_{x_k}B_{ij} [\vec{Y}^n_{M},\tau_n])\Big\}  (\Delta W^n_i \Delta W_j^n - \tau \delta_{ij})
\end{align}
In this equation the extra term with respect to \eqref{sol_BSDE_discret} involves 
$$
\{\cdots \}=\nabla_{\vec x} g[\vec{Y}^n_{M},\tau_n]^T B[\vec{Y}^n_{M},\tau_n] \nabla_{\vec x}\cdot B[\vec{Y}^n_M, \tau_n]
$$ 
(where the divergence is an inner product with the line index of $B$). Therefore, if we are willing to use the explicit expression of $\nabla\cdot B$, we only have to evaluate $\nabla g^T B$ and can revert to the architecture of Fig. \ref{fig:DNN_EM}, with no need of new sub-networks.

\begin{figure}
    \centering
    \includegraphics[width=0.9\columnwidth]{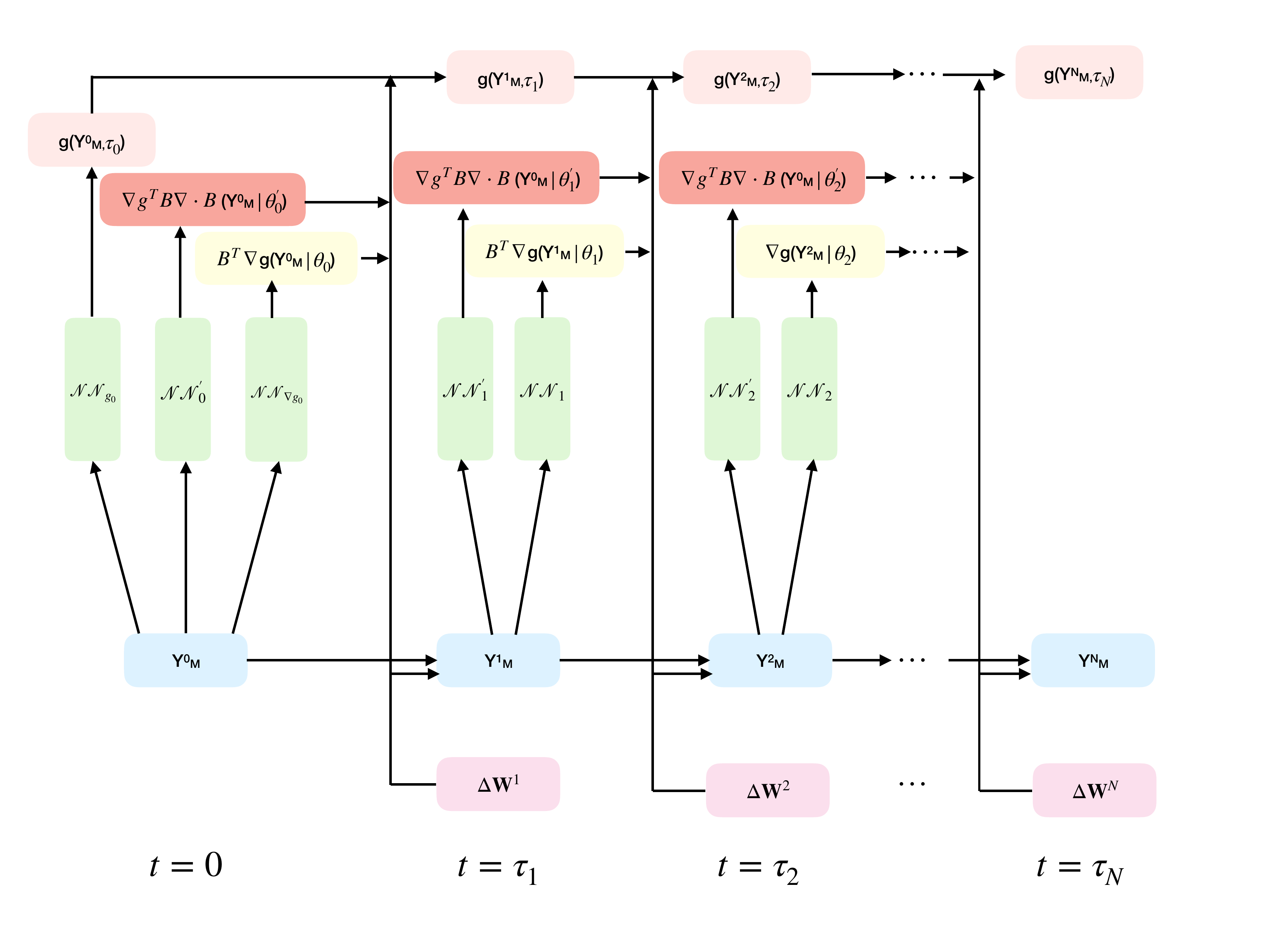}
    \caption{\small $\mathcal{DNN}_{M}$: the network obtained from the Milstein scheme. Time flows horizontally, from left to right, while vertically, the subnetworks output approximations of $g$, $\nabla g^T B$ and $\nabla g^T B \nabla\cdot B$ at each time step. $\mathcal{DNN}_{M}$ takes as {\it input} stochastic trajectories $\{\vec{Y}^{n}_{M}, \Delta{\vec W}^n\}_{0\leq n \leq N}$ and {\it outputs} the value of $\hat{g}(\{\vec{Y}^{n}_{M}, \Delta{\vec W}^n\}_{0\leq n \leq N}\vert \vec{\Theta})$. Specifically, at time $t=0$ the data are processed by two separate sub-neural networks $\mathcal{NN}_{g_{0}}$, $\mathcal{NN}_{\nabla g_{0}}$ and by
    $\mathcal{N}\mathcal{N}_0^\prime$. Their output allows to compute the approximate value of $g(\vec{Y}^{1}_{M},\tau_1)$ using the recursive equation \eqref{sol_BSDE_discret_MILSTEIN}.
     At this point time flows and at $t=\tau_1$ and $(\vec{Y}^{1}_{M}, \Delta{\vec W}^1)$ are given to the sub-networks $\mathcal{NN}_1$ and $\mathcal{NN}_1^\prime$. These returns $\nabla g^T B[\vec{Y}^1_{M}\vert \vec{\theta_1}]$ and
     $\nabla g^T B \nabla\cdot B[\vec{Y}^1_{M}\vert \vec{\theta}_1^\prime]$
     used to get the approximate value of $g(\vec{Y}^{2}_{M},\tau_2)$. This procedure continues till time $\tau_N$ where the approximate value for $g(\vec{Y}^{N}_{M},\tau_N)$ is obtained (see \eqref{estimate}). The internal architecture for 
     $\mathcal{NN}_{g_{0}}$, $\mathcal{NN}_{\nabla g_{0}}$, $\mathcal{NN}_n$ is the same as for the EM case. For the extra networks 
     $\mathcal{NN}_n^\prime$ have one input layer with $d$ neurons, two internal layers with $d+10$ layers, and one output layer
     with $d^2$ neurons to accommodate the $d\times d$ matrix output $\nabla g^T B \nabla\cdot B$.}
    \label{fig:DNN_M}
\end{figure}

However our real interest is in {\it not using} the explicit expression of the diffusion tensor $B$ in the construction of the neural network. There are at least two reasons: i) it is conceptually interesting to avoid the use of the drift term $A$ and diffusion tensor $B$ in the neural architecture; ii) one may imagine applications where good quality input samples are known (the discrete stochastic trajectory) but $A$ and $B$ are not known. Therefore we construct an architecture, which uses two type of sub-networks. One type to approximate 
the vector $\nabla_{\vec x} g[\vec{Y}^n_{M},\tau_n]^T B[\vec{Y}^n_{M},\tau_n]$ and a second type to approximate the matrix $\{\cdots\} =
\nabla_{\vec x} g[\vec{Y}^n_{M},\tau_n]^T B[\vec{Y}^n_{M},\tau_n] \nabla_{\vec x}\cdot B[\vec{Y}^n_{M},\tau_n]$. The later task is achieved by extra sub-networks $\mathcal{NN}_{n}^\prime$ depending on a set of parameter $\vec{\theta}^{'}_n$ for $n=1, \cdots, N-1$ which output an approximate value 
\begin{align}
    (\nabla g^T B \nabla\cdot B)[\vec{Y}^n_M\vert \vec{\theta}_n^\prime] \approx \nabla_{\vec x} g[\vec{Y}^n_{M},\tau_n]^T B[\vec{Y}^n_{M},\tau_n] \nabla_{\vec x}\cdot B[\vec{Y}^n_{M},\tau_n]
\end{align}
The final architecture is a simple modification of $\mathcal{D}\mathcal{N}\mathcal{N}_{EM}$ to include these extra sub-networks at each time step, and is depicted on Fig. \ref{fig:DNN_M}. The full network $\mathcal{D}\mathcal{N}\mathcal{N}_{M}$ inputs the trajectory $\{\vec{Y}^{n}_{M}, \Delta{\vec W}^n\}_{0\leq n \leq N}$ obtained from Milstein discretization scheme \eqref{discretiz_M}
and outputs the estimate $\hat g(\vec{Y}^{n}_{M}, \Delta{\vec W}^n\vert \vec{\Theta})$ according to the process of Fig. \ref{fig:DNN_M}. Again, the set of all parameters is learned by minimizing the loss \eqref{loss_function}.

As a sanity check the simulations in Sect. 5 confirm that the results of the two methods (i.e., using the expression of $\nabla\cdot B$ versus not using it) lead to indistinguishable results.

\subsection{Architecture based on Leimkuhler-Matthews discretization}\label{ALMD}
\begin{figure}
    \centering
    \includegraphics[width=0.9\columnwidth]{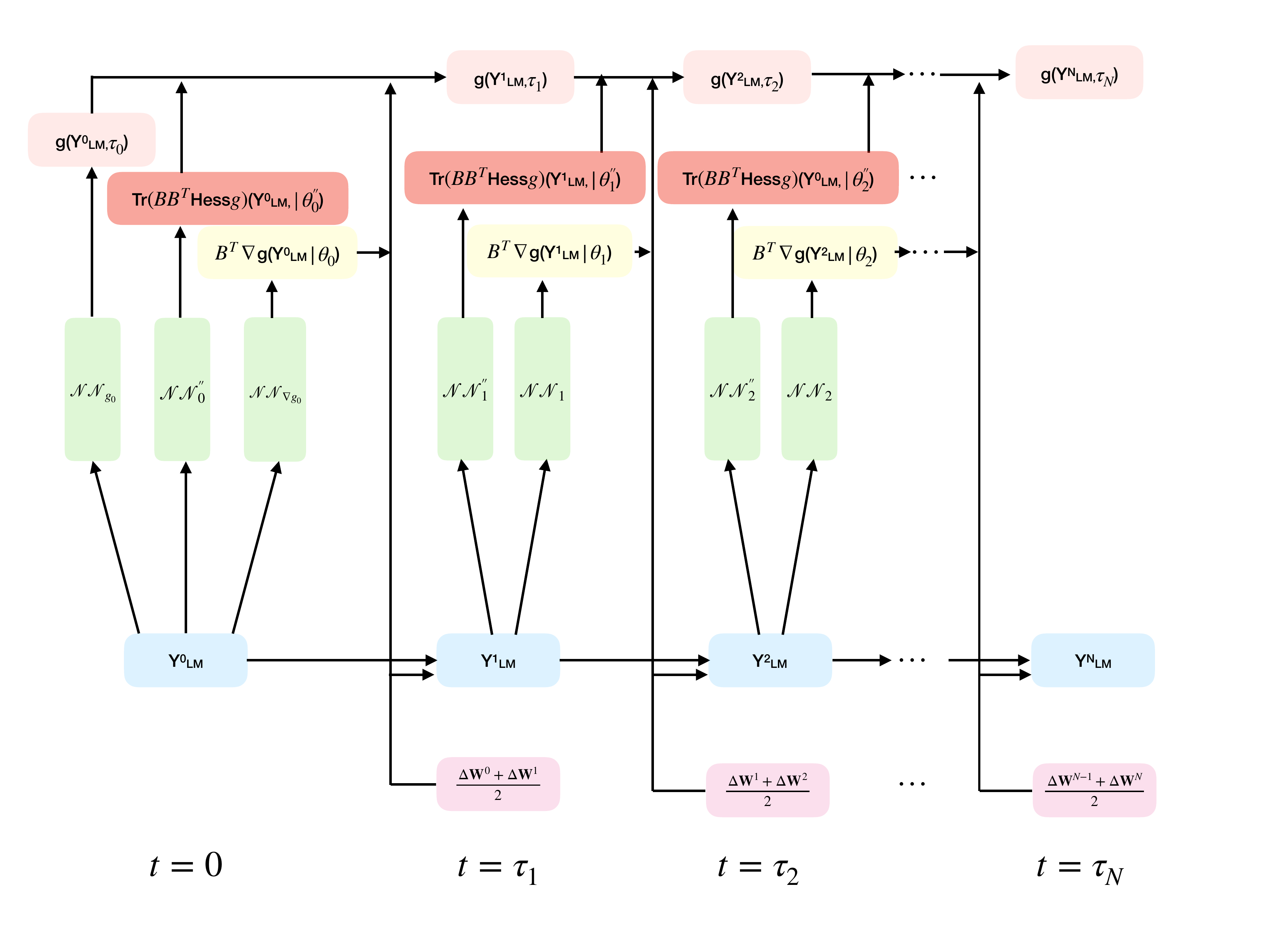}
\caption{\small $\mathcal{DNN}_{LM}$: the network obtained from the Leimkhuler-Matthews scheme. Time flows horizontally while vertically the sub-networks output approximations of $g$, $\nabla g^T B$ and ${\rm Tr} (B B^T {\rm Hess}\, g)$ at each time step. $\mathcal{DNN}_{LM}$ takes as {\it input} stochastic trajectories $\{\vec{Y}^{n}_{LM}, \Delta{\vec W}^n\}_{0\leq n \leq N}$ and {\it outputs} the value of $\hat{g}(\{\vec{Y}^{n}_{LM}, \Delta{\vec W}^n\}_{0\leq n \leq N}\vert \vec{\Theta})$. Specifically, at time $t=0$ the data are processed by three separate sub-neural networks $\mathcal{NN}_{g_{0}}$, $\mathcal{NN}_{\nabla g_{0}}$ and 
    $\mathcal{N}\mathcal{N}_0^{\prime\prime}$. Their output allows to compute the approximate value of $g(\vec{Y}^{1}_{LM},\tau_1)$ using the recursive equation \eqref{sol_BSDE_discret_LM}.
     At this point time flows and at $t=\tau_1$ and $(\vec{Y}^{1}_{LM}, \Delta{\vec W}^1)$ are given to the sub-networks $\mathcal{NN}_1$ and $\mathcal{NN}_1^{\prime\prime}$. These returns $(\nabla g^T B)[\vec{Y}^{1}_{LM}, \tau^1\vert \vec{\theta_1}]$ and 
     ${\rm Tr}(B B^T {\rm Hess}) \,g[\vec{Y}^1_{LM}\vert \vec{\theta}_1^{\prime\prime}]$
     used to get the approximate value of $g(\vec{Y}^{2}_{LM},\tau_2)$. This procedure continues till time $\tau_N$ where the approximate value for $g(\vec{Y}^{N}_{LM},\tau_N)$ is obtained. The internal architecture for 
     $\mathcal{NN}_{g_{0}}$, $\mathcal{NN}_{\nabla g_{0}}$, $\mathcal{NN}_n$ is the same as for the EM case. For the extra network 
     $\mathcal{NN}_n^{\prime\prime}$ we have one input layer with $d$ neurons, two internal layers with $d+10$ layers, and one output neuron. The overall network is trained using SGD to minimize the loss function \eqref{loss_function}.}
    \label{fig:DNN_LM}
\end{figure}

Using the discretization \eqref{discretiz_LM} equation \eqref{sol} becomes (see Appendix-\ref{Appendix-A}) 
\begin{align}\label{sol_BSDE_discret_LM}
    g[\vec{Y}^n_{LM}, & \tau_{n+1}] = g[\vec{Y}^n_{LM},\tau_n]
    \nonumber \\ &
     -f(\tau_n, \vec{Y}^n_{LM}, g[\vec{Y}^n_{LM},\tau_n], B[\vec{Y}^n_{LM},\tau_n]^T \nabla_{\vec x} g[\vec{Y}^n_{LM},\tau_n])\tau
    \nonumber \\ &
    +\frac{1}{2}\nabla g[\vec{Y}^n_{LM},\tau_n]^T B[\vec{Y}^n_{LM},\tau_n] (\Delta\vec{W}^n+\Delta \vec{W}^{n+1})
    \nonumber \\ & 
    - \frac{1}{4}\text{Tr}(B B^T[\vec{Y}^n_{LM}(\tau_{n}),\tau_n] \text{Hess}_{\vec{x}}g[\vec{Y}^n_{LM}(\tau_{n}),\tau_n])\tau 
\end{align}
where, for a sake of simplicity in exposition, we assume that $BB^T=B^TB$.
This time, besides the usual sub-networks for the terms involving $\nabla_{\vec{x}} g$,  we also need extra ones to approximate the term involving the hessian ${\rm Hess}_{\vec x} g$.
The architecture of the full network $\mathcal{D}\mathcal{N}\mathcal{N}_{LM}$ is shown in Fig. \ref{fig:DNN_LM}. With respect to the one of Fig. \ref{fig:DNN_EM} it uses the new sub-networks $\mathcal{N}\mathcal{N}^{''}_n$ that outputs 
\begin{align}\label{hess}
{\rm Tr}(BB^T {\rm Hess} \,g)[\vec{Y}^n_{LM}\vert \vec{\theta}^{''}_n]\approx \text{Tr}(B B^T[\vec{Y}^n_{LM}(\tau_{n}),\tau_n] \text{Hess}_{\vec{x}}g[\vec{Y}^n_{LM}(\tau_{n}),\tau_n]).
\end{align}
where the parameters 
$\vec{\theta}_n^{''}$ have to be learned. These are constituted of one input layer with $d$ neurons, two internal layers with $d+10$ layers, and one output neuron since the output is a scalar. 

\section{Numerical Results}
\label{numerical_result}
In this section, we present the results we have obtained through the use of the three different discretizations. All codes have been developed in Python using numpy, scipy, and Tensorflow libraries. The codes were run on a MacBook Pro 2.4GHz Intel Core i5 with 16GB of RAM. They can be downloaded from \cite{Marino2019GITHUB}. We use the ADAM optimizer \cite{kingma2014adam} for learning the parameters in our numerical examples. 
We have chosen as activation function of neural networks the rectifier, i.e.  $\text{ReLU}(x)=\max(0,x)$.
The algorithm has been successfully tested, for each discretization, on a heat equation in one dimension. In all cases, the algorithm reproduced the expected solution to the requested domain.

We compare the performance of the different networks on two examples examined in \cite{han2018solving} since these have been precisely benchmarked. Moreover, we test the performance of the algorithm by doing an analysis on a non-linear diffusion equation, where an exact solution is known.  

\subsection{Black-Scholes equation}\label{sec:BS}
\begin{figure}
    \centering
    \includegraphics[width=0.9\columnwidth]{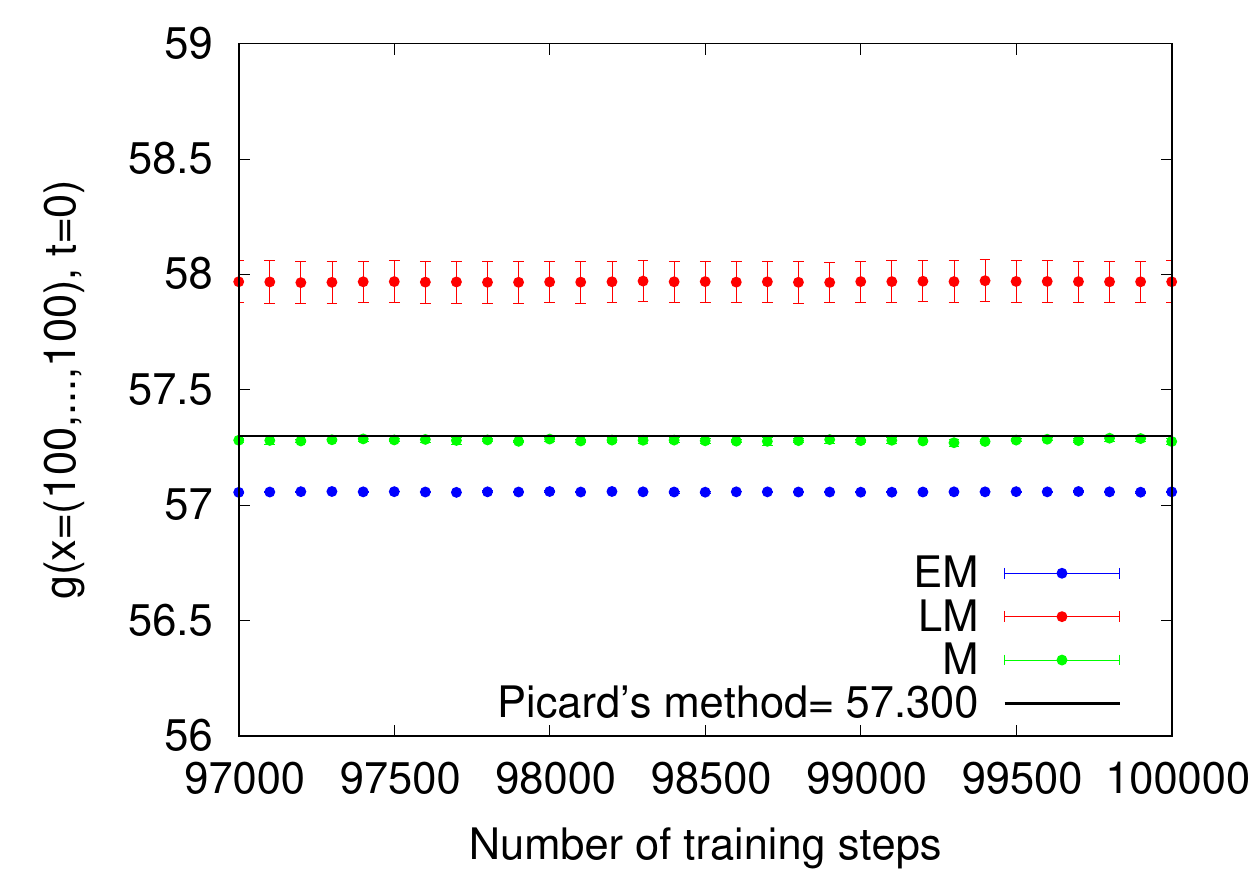}
    \caption{\small The figure shows the value of $g(\vec{x}=(100, \dots, 100), 0)$ returned by the neural networks $\mathcal{DNN}_{EM}$ (blue points), $\mathcal{DNN}_{M}$ (green points) and $\mathcal{DNN}_{LM}$ (red points)  as function of the number of training steps done for learning the parameters. Each point is the mean over five independent runs. The number of equidistant time steps was fixed at $40$, i.e. for each neural networks $\mathcal{DNN}_{EM}$, $\mathcal{DNN}_{M}$ and $\mathcal{DNN}_{LM}$ the value of $N$ is fixed to $40$. All the parameters were initialized randomly uniformly between $[-1,1]$. The total number of training steps, i.e. $t_s$, was fixed at $10^5$ and a dynamic learning rate that goes from $\eta=10$ to $\eta=0.1$ in the following manner: -$\mathcal{DNN}_{EM}$ $\eta_{EM}=10,\,\,\, t_s\in [0, 10^3]$; $\eta_{EM}=5,\,\,\, t_s\in (10^3, 7.5\,10^3]$; $\eta_{EM}=4,\,\,\, t_s \in (7.5\,10^3, 10^4]$; $\eta_{EM}=3,\,\,\, t_s\in (10^4, 2\, 10^4]$; $\eta_{EM}=2,\,\,\, t_s\in (2\, 10^4, 4\, 10^4]$; $\eta_{EM}=1,\,\,\, t_s \in (4\, 10^4, 6.5\, 10^4]$; and $\eta_{EM}=0.1,\,\,\, t_s \in (6.5\, 10^4, 10^5]$; 
    -$\mathcal{DNN}_{M}$ $\eta_{M}=10,\,\,\, t_s\in [0, 10^3]$; $\eta_{M}=5,\,\,\, t_s\in (10^3, 7.5\,10^3]$; $\eta_{M}=4,\,\,\, t_s \in (7.5\,10^3, 10^4]$; $\eta_{M}=3,\,\,\, t_s\in (10^4, 2\, 10^4]$; $\eta_{M}=2,\,\,\, t_s\in (2\, 10^4, 4\, 10^4]$; $\eta_{M}=1,\,\,\, t_s \in (4\, 10^4, 6.5\, 10^4]$; and $\eta_{M}=0.1,\,\,\, t_s \in (6.5\, 10^4, 10^5]$;
    -$\mathcal{DNN}_{LM}$ $\eta_{LM}=10,\,\,\, t_s\in [0, 6.5\,10^4]$; $\eta_{M}=0.1,\,\,\, t_s\in (6.5\,10^4, 10^5]$. The black horizontal line identifies the multilevel Picard method value.}
    \label{fig:BS-comp}
\end{figure}

A prominent example discipline where high-dimensional backward non-linear Kolmogorov equations are used is  financial mathematics. The Black-Scholes partial differential equation \cite{black1973pricing, GVK563580607, gardiner1985handbook} which governs the evolution of the price of options and derivatives, is difficult to solve in high-dimensions, e.g., for a portfolio consisting of $100$ different options, by using classical methods. Numerically, indeed, the time needed in a Monte-Carlo simulation \cite{hammersley2013monte}, to approximate a solution {\it in a whole domain} grows exponentially with the dimension of the portofolio \cite{berner2018analysis}.


Consider the fair price of a European claim based on $100$ underlying assets conditional on no default having occurred yet. When the default of the claim's issuer occurs, the claim's holder receives only a fraction $\delta \in [0,1)$ of the current value. The possible default is modelled by the first jump time of a Poisson process with intensity $Q$, a decreasing function of the current value. In other words, the default becomes more likely when the claim value is low. Mathematically, this model is described by 
\begin{equation}
\label{BS-eq}
\begin{split}
    &\frac{\partial g}{\partial t}(\vec{x},t) + \overline{\mu}\vec{x}\cdot \nabla_{\vec x} g(\vec{x},t) +\frac{\overline{\sigma}^2}{2}\sum^{d}_{i=1} |x_i|^2 \frac{\partial^2 g}{\partial x_{i}^{2}}(\vec{x},t)+ \\
    &-(1-\delta)Q(g(\vec{x},t))g(\vec{x},t)-R\,g(\vec{x},t)=0
\end{split}
\end{equation}
where $R$ is the interest rate of the risk free asset, and $Q(y)$ is a piece-wise linear function given by:
\begin{equation}
\label{Q}
\begin{split}
    &Q(y)= \text{ReLU}\left(\text{ReLU}(y-v_h)\frac{\gamma_h-\gamma_l}{v_h-v_l} + \gamma_h -\gamma_l \right) + \gamma_l.
\end{split}
\end{equation}
Numerically we have chosen the values as: $\delta=2/3,\, R=0.02, \, \overline{\mu}=0.02, \, \overline{\sigma}=0.2,\, v_h=50,\, v_l=70, \, \gamma_h=0.2, \, \gamma_l=0.02$. We fix the dimension to $d=100$, and, therefore, we impose the equation to live on the space $[0,T]\times \mathbb{R}^{100}$, with $T=1$, and we fix the terminal condition be $g(\vec{x}, T)=\phi(\vec{x})=\min\{{x_1,\dots, x_{100}\}}$. Because the exact solution of equation  \ref{BS-eq} is not known, we confine ourselves on a single point $\vec{x}=(100, \dots, 100)$, where analysis performed with the multilevel Picard method \cite{hutzenthaler2017multilevel, weinan2019multilevel} indicates that the value of $g_{\rm Picard}(\vec{x}=(100, \dots, 100), 0)\approx 57.300$ \cite{han2018solving} (black line in Fig. \ref{fig:BS-comp}). 

Fig. \ref{fig:BS-comp} shows the test for the three different networks considered here. Each point is the average over 5 samples. Blue points identify the output of $\mathcal{D}\mathcal{N}\mathcal{N}_{\rm EM}$ and reach the approximate value $g_{\rm EM}(\vec{x}=(100, \dots, 100), 0)=57.059 \pm 0.003$ in agreement with \cite{han2018solving}. Red points show that the value obtained from $\mathcal{D}\mathcal{N}\mathcal{N}_{\rm LM}$,  $g_{\rm LM}(\vec{x}=(100, \dots, 100), 0)=57.968 \pm 0.091$, is substantially larger with respect to the Euler-Maruyama and the multilevel Picard method. The green points represent the value obtained by $\mathcal{D}\mathcal{N}\mathcal{N}_{\rm M}$ $g_{M}(\vec{x}=(100, \dots, 100), 0)=57.276 \pm 0.012$, which agrees with Picard's method value within $3$ standard deviation of the mean. 

We see that the network based on Milstein's discretization lead to improved accuracy with a computational complexity similar to the one of the Euler-Maruyama discretization. The network based on the Leimkuhler-Mattews discretization on the other hand do not show good performance. A priori this could come from the use of the extra sub-network to estimate the terms containing ${\rm Tr} BB^T{\rm Hess}_{\vec x} g$, or could also come from the fact that the time step $\tau$ has to be taken quite smaller.

\begin{figure}
    \centering
    \includegraphics[width=0.9\columnwidth]{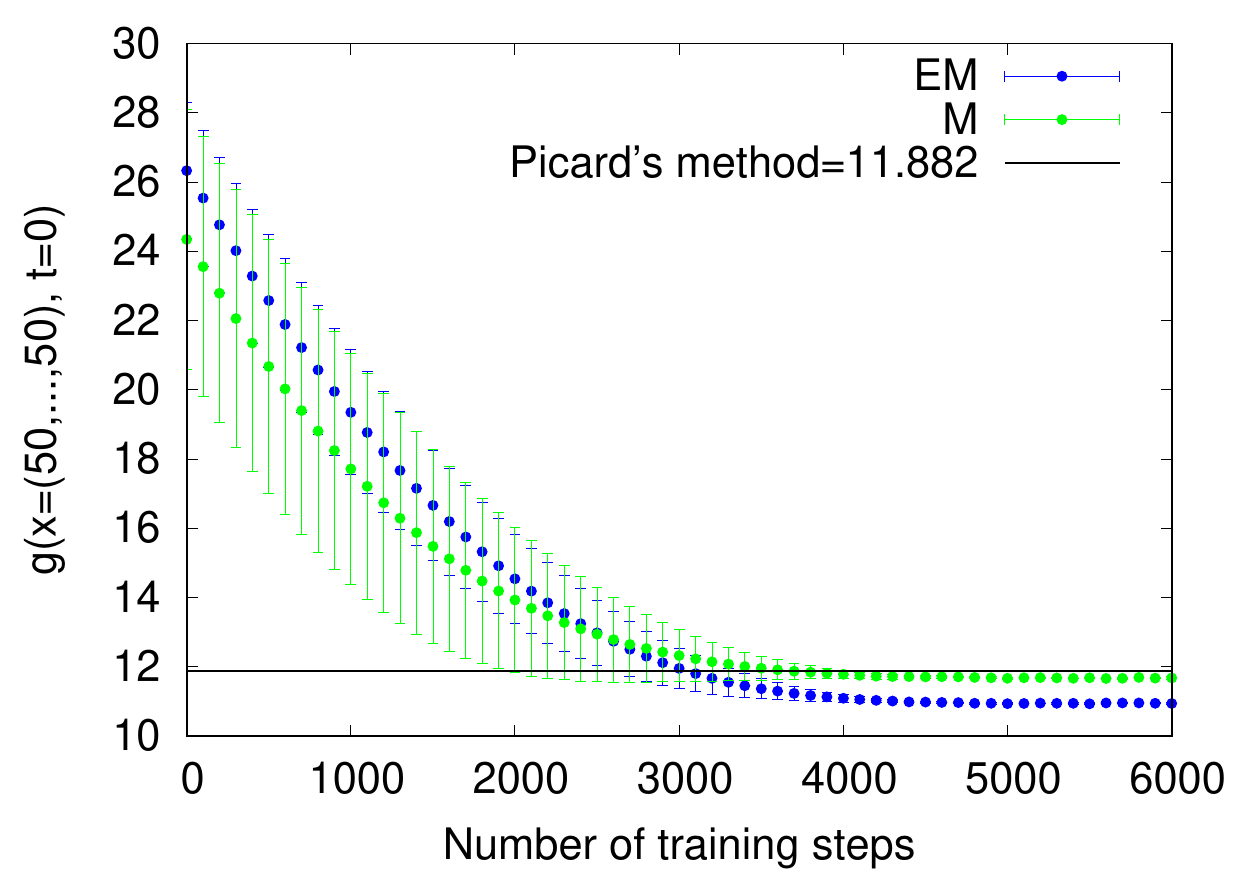}
    \caption{\small The figure shows the value of $g(\vec{x}=(50, \dots, 50), 0)$ returned by the neural networks $\mathcal{DNN}_{EM}$ (blue points), $\mathcal{DNN}_{M}$ (green points) as function of the number of training steps done for learning the parameters. Each point is the mean over five independent runs. The number of equidistant time steps was fixed at $40$, i.e. for each neural networks $\mathcal{DNN}_{EM}$ and $\mathcal{DNN}_{M}$  the value of $N$ is fixed to $40$. All the parameters were initialized randomly uniformly between $[-1,1]$. The total number of training steps, i.e. $t_s$, was fixed to $6000$ and the learning rate  to $\eta=0.008$. The black horizontal line identifies the multilevel Picard method value.}
    \label{fig:BSNL-comp}
\end{figure}

Observing a better performance for the Milstein scheme, we decided to re-test it on another Black-Scholes equation. This time we made a comparison only with the Euer-Maruyama scheme. The equation is similar to the one in (\ref{BS-eq}). It can be read:
\begin{equation}
\label{BSNL-eq}
\begin{split}
    &\frac{\partial g}{\partial t}(\vec{x},t) + \vec{x}\cdot \nabla_{\vec x} g(\vec{x},t) +\frac{1}{2}\sum^{d}_{i=1} |x_i|^2 \frac{\partial^2 g}{\partial x_{i}^{2}}(\vec{x},t)+ \exp ( -g(\vec{x},t))=0.
\end{split}
\end{equation}
We fix the terminal condition to be $g(\vec{x}, T=0.5)=\phi(\vec{x})=\min\{{x_1,\dots, x_{100}\}}$. Again,  because the exact solution of equation  \ref{BSNL-eq} is not known, we confine ourselves on a single point $\vec{x}=(50, \dots, 50)$, where analysis performed with the multilevel Picard method \cite{hutzenthaler2017multilevel, weinan2019multilevel} indicates that the value of $g_{\rm Picard}(\vec{x}=(50, \dots, 50), 0)\approx 11.882$ (black line in Fig. \ref{fig:BSNL-comp}). 
In Fig. \ref{fig:BSNL-comp} we present the comparison between the output of $\mathcal{D}\mathcal{N}\mathcal{N}_{\rm EM}$ (blue points) and $\mathcal{D}\mathcal{N}\mathcal{N}_{\rm M}$ (green points). 
Each point is the average over 5 samples. Blue points identify the output of $\mathcal{D}\mathcal{N}\mathcal{N}_{\rm EM}$ and reach the approximate value $g_{\rm EM}(\vec{x}=(50, \dots, 50), 0)=10.942 \pm 0.027$. The green points represent the value obtained by $\mathcal{D}\mathcal{N}\mathcal{N}_{\rm M}$ $g_{M}(\vec{x}=(50, \dots, 50), 0)=11.678 \pm 0.032$, which agrees with Picard's method value. We observe again that the network based on Milstein's discretization lead to improved accuracy with a computational complexity similar to the one of the Euler-Maruyama discretization.

In the following paragraph we test this issue by considering a simpler equation which has a constant diffusion matrix.

\subsection{Allen-Cahn equation}\label{sec:AC}

\label{ACeq} 
\begin{figure}
    \centering
    \includegraphics[width=0.45\columnwidth]{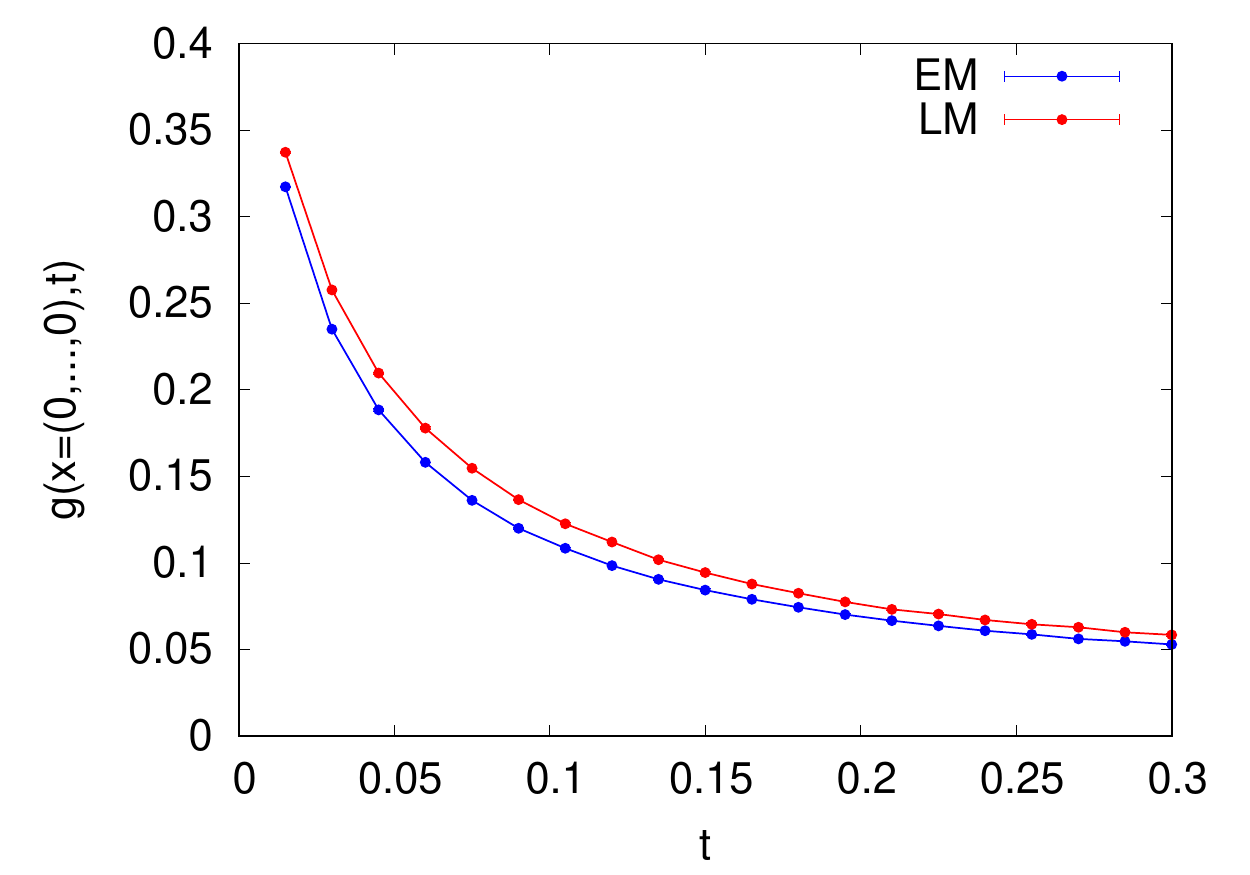}
    \includegraphics[width=0.45\columnwidth]{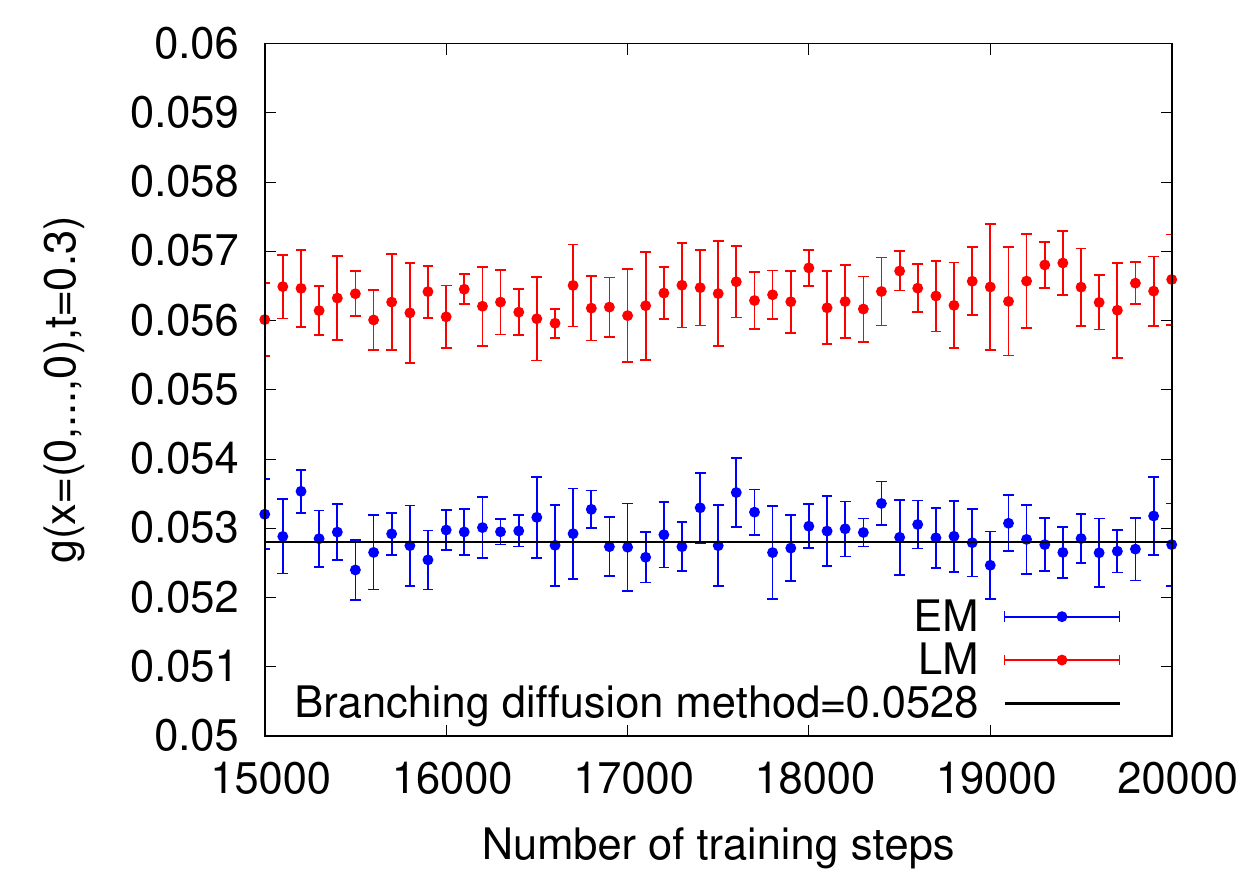}
    \caption{\small Comparison of the ouputs of $\mathcal{D}\mathcal{N}\mathcal{N}_{\rm EM}$ (blue points) and  
    $\mathcal{D}\mathcal{N}\mathcal{N}_{\rm LM}$ (red points) for the Allen-Cahn equation \eqref{AC} with $\epsilon=1$ in d=100 dimensions. \textbf{Left}: the figure shows the temporal evolution of the solution of the equation \ref{AC} for the two methods at a fixed point $\vec{x}=(0,...0)$. \textbf{Right}: the figure shows in detail the discrepancy between the LM network and EM network in reaching a precise value for $g(\vec{x}=(0,....,0))$ at time $t=0.3$. The EM result agrees in reaching the predicted value obtained by the branching diffusion method. Once again, the Leimkuhler-Matthews  method seems to overestimate the value of the solution.} 
    \label{fig:ACcomp}
\end{figure}

\begin{figure}
    \centering
    \includegraphics[width=0.45\columnwidth]{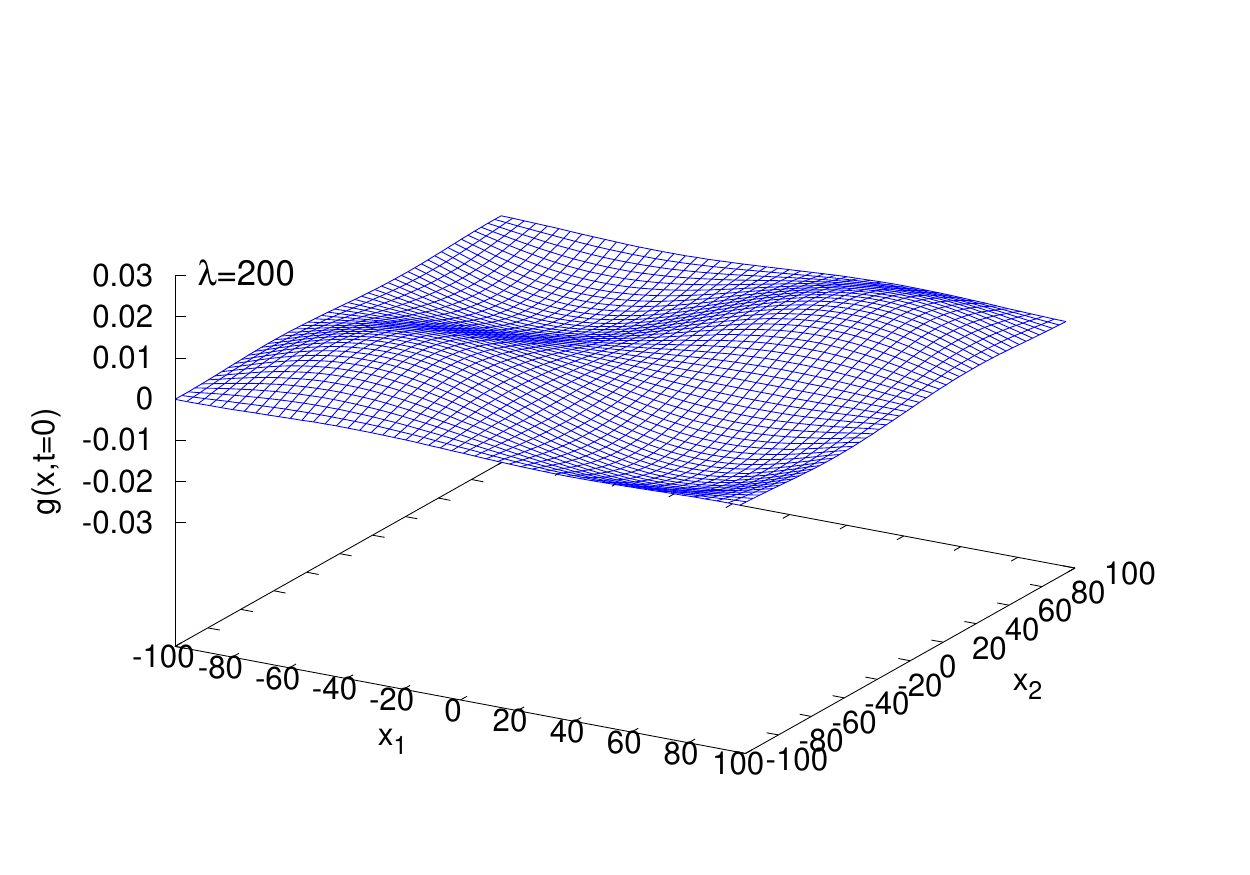}
    \includegraphics[width=0.45\columnwidth]{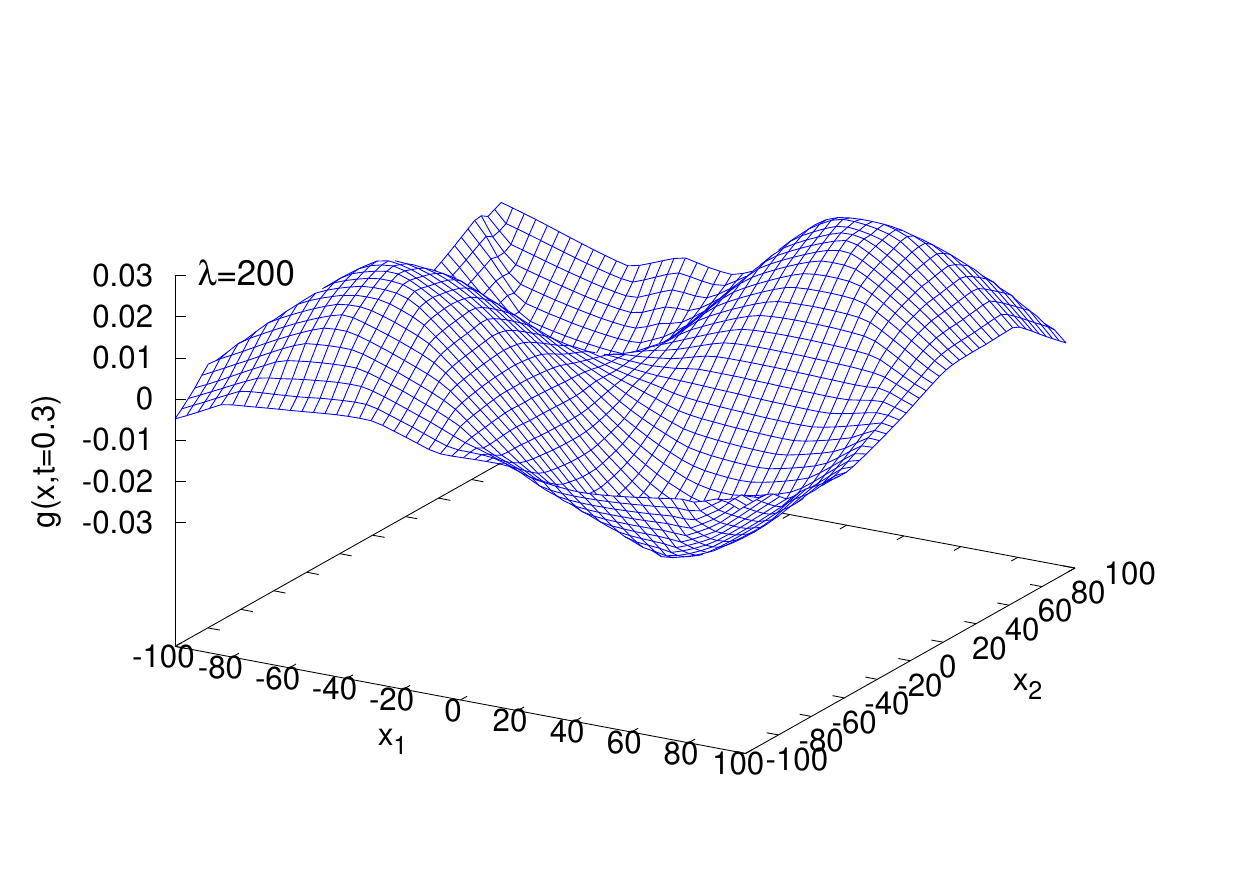}
    \caption{\small Allen-Cahn equation with $\epsilon=0.5$ in the plane $[-100,100]^2$. The figure shows the growth of an initial condition (left) $\Phi(x_1, x_2) = \delta_0 \sin(\frac{2\pi x_1}{\lambda}) \sin(\frac{2\pi x_2}{\lambda})$ where $\lambda=200$ and $\delta_0=0.01$. The right figure shows the approximate solution  $g_0(x_1, x_2\vert\vec{\theta_0})$ given by the deep networks (computed on $10^4$ points) at time $t=0.3$. As expected, since $\lambda\gg\lambda_c \approx 4.44$ the pattern grows (right), and is attracted by the two minima of the potential. When $\lambda\ll\lambda_c$ one checks that the initial condition is attracted to zero over the whole domain.} 
    \label{fig:ACfull}
\end{figure}

The Allen–Cahn equation is a reaction–diffusion equation that arises in Ginzburg-Landau theories of order–disorder transitions and pattern formation. The order parameter $g(\vec{x},t) \in \mathbb{R}$ evolves according to the forward semi-linear parabolic PDE
\begin{equation}
\label{AC}
\frac{\partial}{\partial t} g(\vec{x}, t) = \Delta g(\vec{x}, t) + \epsilon^{-2}( g(\vec{x}, t) - g(\vec{x}, t)^3)
\end{equation}
with {\it initial} condition $g(\vec{x}, 0) = \phi(\vec{x})$. Here $\Delta$ is the $d$-dimensional Laplacian.
The equation describes the competition between a tendency towards phase separation towards two attractive minima of a potential 
$\propto (1-g^2)^2$ and spreading due to diffusion arising from a surface tension $\propto\vert\nabla g\vert^{2}$. When $\epsilon\to 0$ the equation describes motion by curvature of thin interfaces. On the other hand when $\epsilon\to +\infty$ we get the heat equation which erases any order and pattern (see Appendix-\ref{Appendix-B} for a the analysis of the heat equation).

The transformation $t\to T -t$ maps the Allen-Cahn equation to a backward Kolmogorov non-linear equation of the form \eqref{Kolmogoroveq} with 
zero drift $A=0$, constant diagonal diffusion matrix $B= 2^{1/2}I$ and non linearity $f(g) = \epsilon^{-2}(g - g^3)$. 
It is then straightforward to compare the outputs of $\mathcal{D}\mathcal{N}\mathcal{N}_{\rm EM}$ and $\mathcal{D}\mathcal{N}\mathcal{N}_{LM}$. The result is displayed on Fig. \ref{fig:ACcomp} which shows that even with constant and diagonal diffusion tensor the Leimkuhler-Matthews scheme overestimates the solution of the Allen-Cahn equation. However, we also observe that by taking smaller time step $\tau$ the situation seems to improve slightly. Nevertheless this becomes quickly impractical in terms of memory requirements because the number of sub-networks is proportional to the number of time steps. In this simulation we have taken 
$\phi(\vec{x})=1/(2+0.4||\vec{x}||^2)$, $d=100$, and $\epsilon =1$ as in \cite{han2018solving}.

We stress that the present methodology allows to compute an approximation of the solutions $g_0(\vec{x}\vert \vec{\theta}_0)$ for a {\it whole domain}. For purely illustrative purposes let us consider the Allen-Cahn equation in two dimensions for a periodic initial perturbation. Figure \ref{fig:ACfull} shows the solution of \eqref{AC} for the initial condition $\Phi(\vec{x}) = \delta_0 \sin(\frac{2\pi x_1}{\lambda}) \sin(\frac{2\pi x_2}{\lambda})$ where 
$\delta_0\ll 1$. The figure can be understood by linear stability analysis which yields the critical wavelength $\lambda_c \approx \sqrt{2}\,2\pi \epsilon$ such that smaller wave lengths are damped by the diffusive effects, whereas longer wavelengths grow and are attracted by the two minima of the potential. 

\begin{figure}
    \centering
    \includegraphics[width=0.45\columnwidth]{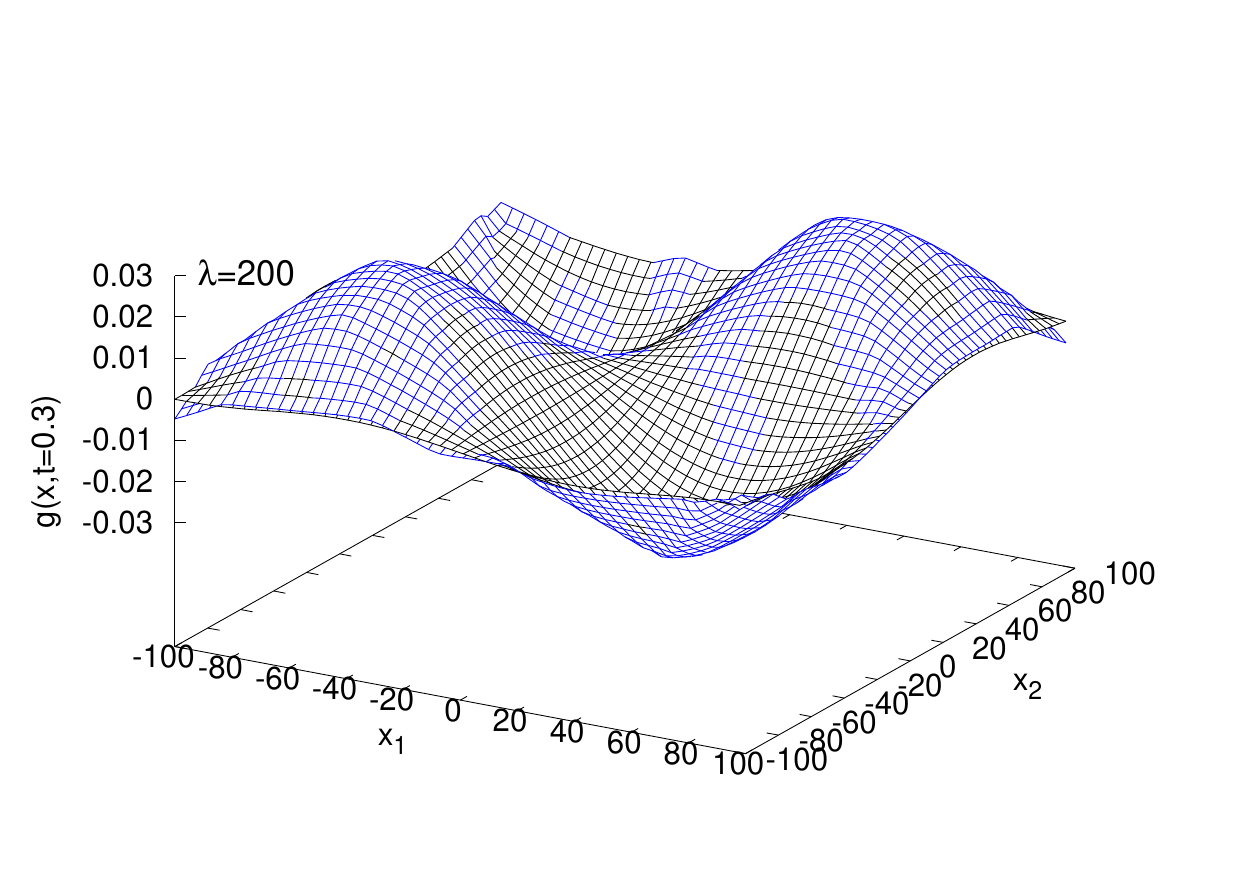}
    \caption{\small 
    Allen-Cahn equation with $\epsilon=0.5$ in the plane $[-100,100]^2$. The figure shows the comparison of two different algorithms for computing the approximate solution  $g_0(x_1, x_2)$ on $10^4$ points at time $t=0.3$. In blue is plotted the solution given by he deep networks algorithm BSDE using EM discretisation, while in black the solution was obtained by using the Numerical Method of Lines, developed on Mathematica. The difference between the two methods is very small, and on average the absolute error is of the order of $0.009\%$.
    }
    \label{fig:ACfull_comparison}
\end{figure}

Moreover, we show in Figure \ref{fig:ACfull_comparison} the comparison of the BSDE using EM discretisation, with a well established algorithm used for solving PDEs with Mathematica Software, i.e., the Numerical Method of Lines (NML). The average absolute error  between the two methods is really small. More precisely, we run the two algorithms, under the same initial condition as in Figure \ref{fig:ACfull},  on the same domain  $[-100,100]^2$, discretised in $M=10^4$ points. Then we defined the average absolute error as:
\begin{equation}
	\frac{1}{M} \sum_{m=1}^M \left| g_{NML}(\vec{x}_m, 0.3) - \mathcal{N}\mathcal{N}(\vec{x}_m, 0.3|\vec{\theta}) \right|,
\end{equation}
and we computed the values over the $M$ points. Surprisingly, the average absolute error is equal to $0.009\%$, showing that the two methods are comparable in approximating the solution of the Allen-Cahn equation in $2d$.

In the experiments above, the reader would expect the comparison also of the Milstein method. However, for the Allen-Cahn equation in (\ref{AC}), it cannot be applied. Indeed, the diffusion coefficient is constant, which implies that the $x$ derivative of $B$, in the third line of (\ref{sol_BSDE_discret_MILSTEIN}), is zero. For filling this gap, we present below a modified version of the Allen-Cahn equation, where the diffusion coefficient is  dependent by $\vec{x}$. More precisely, we modify the equation (\ref{AC}) as:
\begin{equation}
\label{AC_modified}
\frac{\partial}{\partial t} g(\vec{x}, t) = \Delta (\mathbf{D}(x) g(\vec{x}, t)) + \epsilon^{-2}( g(\vec{x}, t) - g(\vec{x}, t)^3)
\end{equation}

where $D(x)_{ii}=(x_{i}^2)$, and $D(x)_{ij}=0$ if $i\neq j$. The initial condition is the same as (\ref{AC}) and $d=100$. We compare the results with the Mutilevel Picard Algorithm  \cite{Weinan2016MultilevelPI, 10.1007/s10915-018-00903-0} .  The parameter $\epsilon$ is set to $1$. The results only for the Euler-Maruyama and Milstein schemes are listed in Tab. $1$. 
On average the Milstein scheme approximates the solution of  eq. (\ref{AC_modified}) better than the Euler-Maruyama scheme. The values in Tab. $1$ were obtained averaging over 5 different experiments.

\begin{table}
\label{tab-dAC100}
\centering
\scalebox{1}{
\begin{tabular}{ |c||c||c||} 
\hline
 Reference Value MLP & Value obtained with EM scheme &  Value obtained with M scheme  \\
\hline
0.55706 & 0.55692 (2)& 0.55705 (1) \\
\hline
\end{tabular}
}
\caption{\small The table shows the average values of $g(\vec{x}=(0.0005, ..., 0.0005), t=0.15)$ over 5 different experiments. The integers in the brackets identify the standard errors of the mean on the last digit. On average Milstein approximation works better than the Euler-Maruyama one.}
\end{table}


\subsection{An exactly solvable non-linear diffusion equation}\label{sec::non-linear-diff-eq}

As illustrated above (for the Allen-Cahn equation) the deep learning methodology studied in this paper allows to compute an approximation of the solutions for a {\it whole domain}. In this paragraph we show that the method does not suffer from the curse of dimensionality in the sense that the number of training trajectories used scales polynomially in the dimension of space. 

First one computes the output function of the neural net $\mathcal{N}\mathcal{N}(\vec{x},0|\vec{\theta})$ (for a backward problem say, for $\vec{x}\in [0,1]^d$ and time $t=0$) by training the deep network with a set of stochastic trajectories starting at $P$ random initial points (the number of trajectories used for each initial point being fixed). Since the given architecture of the neural network involves a polynomial number of edges w.r.t $d$, we are interested to know how $P$ scales with $d$, in order to attain a pre-specified error. For the error measure, we follow \cite{beck2018solving} and consider an {\it average relative error}
\begin{equation}\label{rel_err_eq_DE}
    \langle \epsilon \rangle = \frac{1}{M} \sum_{m=1}^M \left| \frac{g(\vec{x}_m, 0) - \mathcal{N}\mathcal{N}(\vec{x}_m, 0|\vec{\theta})}{g(\vec{x}_m, 0)} \right|,
\end{equation}
where $\vec{x}_m$, $m=1,\cdots, M$ is a large number of points taken uniformly at random in $[0,1]^d$, and $g(\vec{x}_m, 0)$ should be a
{\it reliable} solution of the non-linear Kolmogorov PDE. Note that \eqref{rel_err_eq_DE} is the Monte-Carlo approximation of the corresponding integral over $[0,1]^d$ and thus we expect a fluctuation of roughly $1/\sqrt{M}$ for different choices of random points.

In practice it is quite difficult to compute the average relative error because in high dimensions reliable solutions are costly to compute by classical methods for a large number $M$ of points. Here we use the {\it exactly solvable} non-linear diffusion equation
\begin{equation}\label{diffeq}
    \frac{\partial g(\vec{x}, t)}{\partial t}+ D\Delta_x g(\vec{x}, t)=-\left[2D g(\vec{x}, t)- \frac{1}{d} - D \right]\sum^{d}_{i=1}\frac{\partial g(\vec{x}, t)}{\partial x_i}
\end{equation}
where $D$ is a constant diffusion coefficient (note that the r.h.s of \eqref{diffeq} conserves the total initial mass).  With terminal condition 
\begin{equation}
 g(\vec{x}, T)=\frac{1}{1+e^{-T-\sum^d_{i=1}x_i}}
 \end{equation}
the solution of the backward problem is known $\forall t \in [0,T]$,
namely
\begin{equation}\label{solutiondiffeq}
    g(\vec{x},t)=\frac{1}{1+e^{-t-\sum^d_{i=1}x_i}}.   
\end{equation}

The deep network of Sec.\ref{subsec::EM} (with Euler-Maruyama discretization) is used since the diffusion term is homogeneous. For the number of initial random points in $[0,1]^d$ of stochastic input trajectories we take $P=4, 8,...., 2048, 4096$. Once the network is trained we obtain the 
approximate solution $\mathcal{N}\mathcal{N}(\vec{x},0|\vec{\theta})$.
The relative error is then easily computed from \eqref{rel_err_eq_DE}, \eqref{solutiondiffeq} with $M=10^4$ uniformly random points $\vec{x}_m\in [0, 1]^d$, $m=1, \cdots, M$. 

\begin{figure}
    \centering
    \includegraphics[width=0.45\columnwidth]{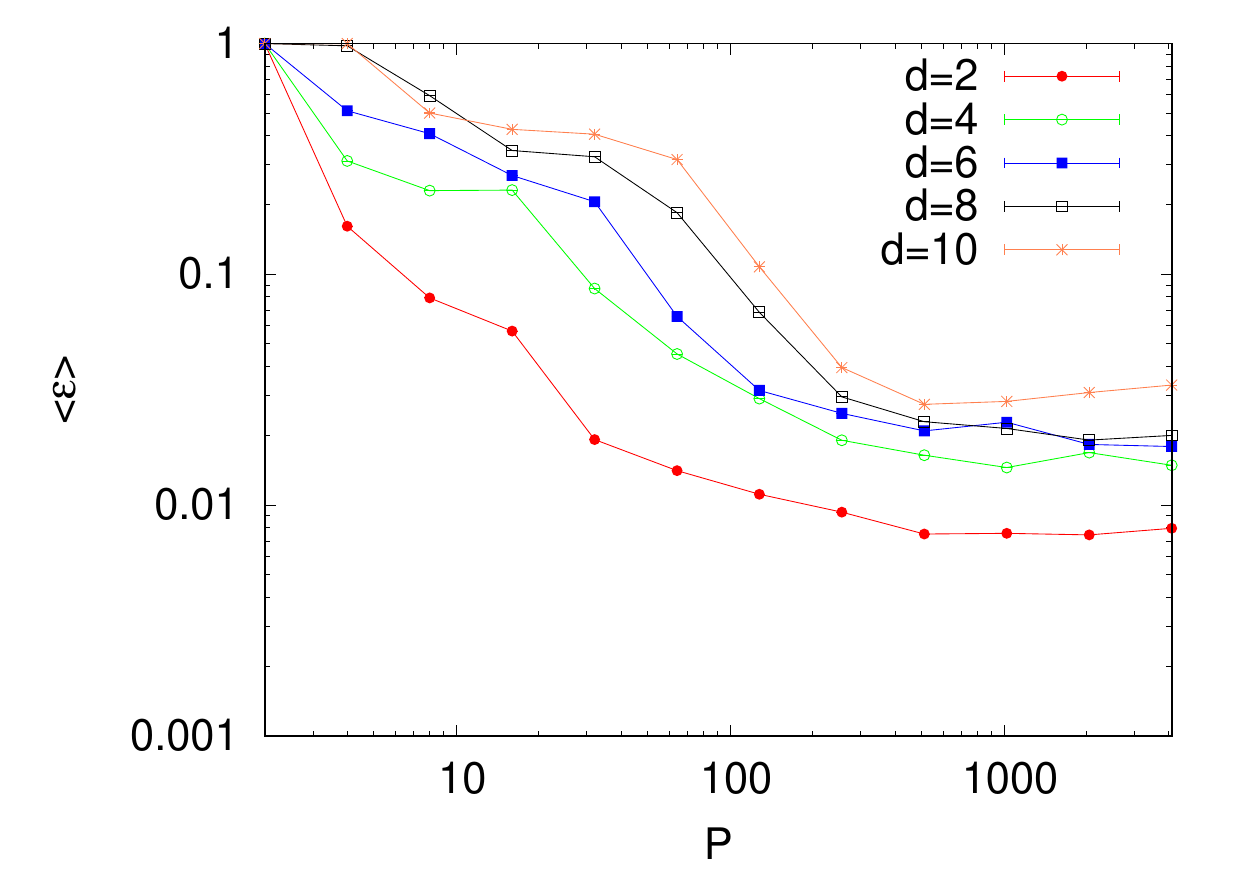}
    \includegraphics[width=0.45\columnwidth]{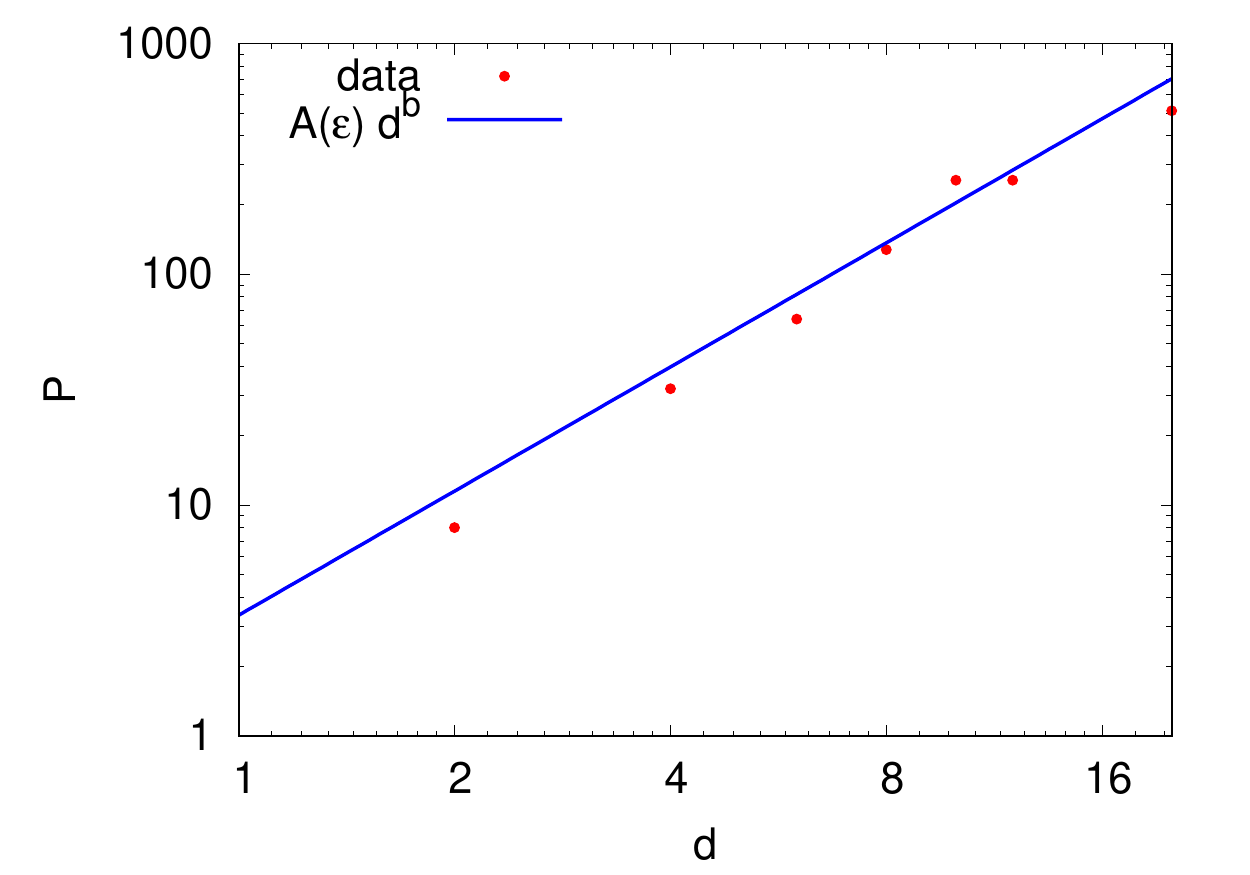}
    \caption{\small \textbf{Left}: Average relative error $\langle \epsilon \rangle$ as a function of number of initial points of stochastic trajectories $P=4,8,16,\cdots, 1024,2048,4096$,
    on a log-log scale i.e. for different values of $d=2,4,6,8,10$.
    \textbf{Right}: Computational complexity, given an average relative error $\langle \epsilon \rangle < 10 \%$, measured by $P$ as a function of $d=2,4,6,8,10,12,18$, on a log-log scale the slope is $b\approx 1.78$. The computational complexity obtained is proportional to $A(\epsilon)d^{1.78} $, with $A(\epsilon) \sim 0.20$. The results are obtained using the algorithm described in Sec.\ref{subsec::EM}, and by fixing $D=\frac{\sigma^2}{2}$, where $\sigma=0.25$, and $T=0.01$.}
    \label{fig-err}
\end{figure}

Fig. \ref{fig-err} shows the behaviour of $\langle\epsilon\rangle$ as a function of $P$, for various values of $d$. For small value of $d$, the averaged relative error needs few points $P$ to approximate well the exact solution. In contrast, as expected, as $d$ grows, the number of points $P$ needed for approximate well the solution becomes larger. Note that the error floor in the relative error equals a few percent which corresponds to the magnitude of fluctuations in 
\eqref{rel_err_eq_DE}.
By fixing a pre-specified relative error of $\langle \epsilon \rangle \approx 0.1$ (well above the error floor) we observe that $P$ scales polynomially with $d$, roughly as $P=O(d^{1.78})$.
In Appendix \ref{Appendix_complexity_Milstein} we present a  similar experiment for the Milstein scheme.
In Appendix \ref{Appendix-time} we also point out that, in practice, the time for approximation also scales polynomially with the dimension.

\subsection{The Hamilton-Jacobi-Bellman equation}\label{sec::HJB-sec}

\begin{figure}
    \centering
    \includegraphics[width=0.9\columnwidth]{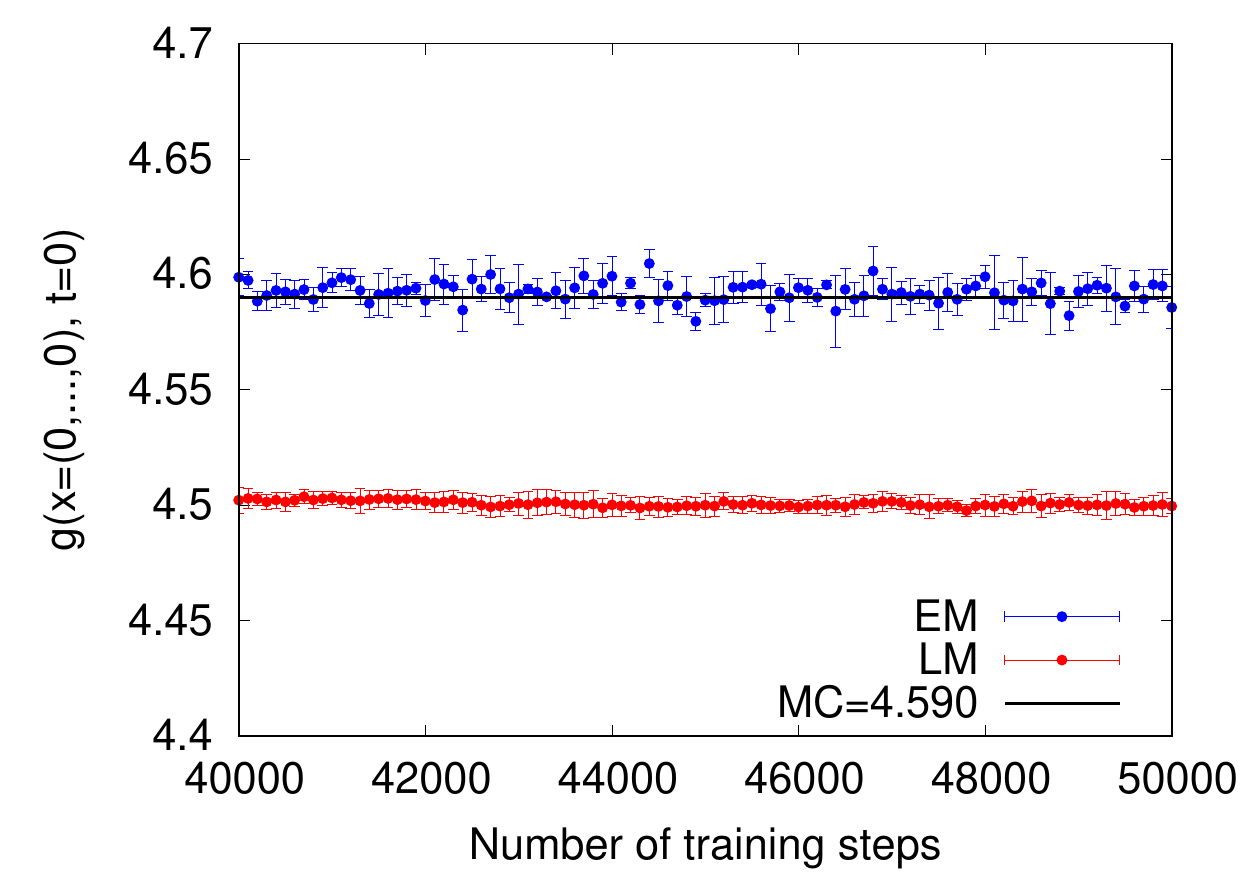}
    \caption{\small The figure shows the comparison between the Euler-Maruyama (blue points) and the Leimkuhler-Mattews (red points) schemes in approximating the solution $g(\vec{x}, t=0)$ of equation \eqref{HJB-eq}, with $\lambda=1$, at the origin of $\mathbb{R}^{100}$, as a function of the training steps. Each point is the average over $5$ experiments. Error bars are standard deviations. The black line is the exact  value obtained by standard Monte Carlo on equation \eqref{HJB-eq_sol}.  }
    \label{fig:HJB-compLMvsEM}
\end{figure}

The Hamilton–Jacobi–Bellman (HJB) equation is used frequently in optimal control theory. It, indeed, gives a necessary and sufficient condition for optimality of a control with respect to a loss function. In general, it is  a nonlinear partial differential equation in the value function, which means its solution is the value function itself. Once the solution is known, it can be used to obtain the optimal control by taking the maximizer/minimizer of the Hamiltonian involved in the HJB equation. The backward HJB equation is defined as:
\begin{equation}
\label{HJB-eq}
    \frac{\partial g}{\partial t}(\vec{x},t) + \Delta g(\vec{x},t)  = \lambda ||\nabla g(\vec{x},t)||^2,
\end{equation}
where $\lambda$ is a positive constant representing the strength of the control. Applying Ito's Lemma to (\ref{HJB-eq}), conditioned to terminal condition $g(\vec{x}, T)=\phi(\vec{x})=\ln((1+||\vec{x}||^2)/2)$ with $\vec{x} \in \mathbb{R}^d$, an exact solution arises naturally. It is:
  \begin{equation}
\label{HJB-eq_sol}
   g(\vec{x}, t)= -\frac{1}{\lambda}\ln\left(\mathbb{E} \left[-\lambda g(\vec{x}+\sqrt{2}\vec{W}_{T-t}) \right] \right),
\end{equation}
where $\vec{W}_t$ is a standard Brownian motion.

For this particular equation we perform an analysis for comparing the Euler-Maruyama and the Leimkuhler-Mattews schemes.
In Fig. \ref{fig:HJB-compLMvsEM} we present the comparison between these two different algorithms. In blue points we present the the average value of the approximate solution $g(\vec{x}, t=0)$ at $\vec{x}=(0,...0) \in \mathbb{R}^{100}$ of equation \eqref{HJB-eq} when $\lambda=1$ over $5$ different experiments as a function of the number of training sample used for learning parameters of the Neural Network obtained by Euler-Maruyama scheme. As the figure shows, these points are able to approximate the exact solution obtained by standard Monte Carlo, presented as black line at $\sim 4.590$. In contrast, red points identify the average value of the approximate solution of equation \eqref{HJB-eq} returned by using the Leimkuhler-Mattews scheme. This scheme, unfortunately, is not able to reproduce correctly, as in the previous numerical examples, the correct value. 
In Appendix \ref{Appendix_HJBNCD} we present a comparison between the Euler-Maruyama and the Milstein scheme on an HJB equation having a non constant diffusion tensor.

\section{Conclusion}

In this article, we investigated properties of a recent deep learning methodology to find solutions of high-dimensional semi-linear PDEs. First, we have shown numerically that the higher order Milstein discretization improves the accuracy without increasing the order of computational complexity with respect to the first order Euler-Maruyama approximation. On the other hand the first order Leimkuhler-Matthews approximation does not perform better than the Euler-Maruyama one, at least for practical values of the time steps. 

Secondly, using an exactly solvable case we have computed the average relative error of the approximate solution and shown that for a given error the number of input stochastic trajectories scales polynomially with the dimension. 

We hope that this study can trigger further systematic explorations of advantages and disadvantages of higher order discretizations in various proposed architectures \cite{hure2019some, DBLP:journals/corr/abs-1908-00412, DBLP:journals/jscic/Chan-Wai-NamMW19, raissi2018forward, beck2019deep}.
Finally let us mention that the representation power of deep learning methods have been rigorously analyzed essentially only for linear equations \cite{berner2018analysis, grohs2018proof},  where the architectures drastically simplify,  and it would be highly desirable to extend such analysis to the non-linear setting (see \cite{Hutzenthaler-2020} for recent results).

\section*{Acknowledgement}
The work of R. M. was supported by Swiss National Foundation grant number 200021E 17554 and now he is supported by the FARE  grant No. R167TEEE7B. N. M. learned about the topic during the conference "Intelligent Machines and Mathematics" in Bologna (January 2019) and acknowledges interesting related discussions with Pierluigi Contucci and Philip Grohs. We also thank Martin Hutzenthaler as well as the referees for constructive comments that improved the paper.

\section*{Data Availability Policy}
The data that support the findings of this study are available from the corresponding author upon reasonable request. 

\section*{Competing interests}
The authors declare no competing financial interests.
\appendices

\section{Derivation of equation (\ref{sol_BSDE_discret_LM})}
\label{Appendix-A}
For completeness we derive equation (\ref{sol_BSDE_discret_LM}) which follows from the Leimkuhler and Matthews  discretization of stochastic trajectories given by \eqref{discretiz_LM}. Note the differences between this equation and \eqref{sol_BSDE_discret} implied by the Euler-Maruyama discretization \eqref{discretiz_EM}.
Using (14) and applying Ito's lemma (i.e., expanding to second order) we get for small enough $\tau$
\begin{equation}
    \label{eq::LM}
    \begin{split}
    & g[\vec{Y}^n_{LM}, \tau_{n+1}]-g[\vec{Y}^n_{LM},\tau_n]= \partial_t g[\vec{Y}^n_{LM},\tau_n]\tau \\
    & + \nabla_{\vec x} g[\vec{Y}^n_{LM},\tau_n]^TA[\vec{Y}^n_{LM},\tau_n]\tau    +\frac{1}{2}\nabla g[\vec{Y}^n_{LM},\tau_n]^T B[\vec{Y}^n_{LM},\tau_n] (\Delta \vec{W}^n+\Delta \vec{W}^{n+1})+\\
    &\frac{1}{8}\,\mathbb{E}\bigl[(\Delta \vec{W}^n+\Delta \vec{W}^{n+1})^T  B^T[\vec{Y}^n_{LM}(\tau_{n}),\tau_n] \text{Hess}_{\vec{x}}g[\vec{Y}^n_{LM}(\tau_{n}),\tau_n]B(\Delta\vec{W}^n+\Delta \vec{W}^{n+1})\bigr].        
    \end{split}
\end{equation}
To evaluate the last term we work in components and use 
$\mathbb{E}[\Delta W_j^n\Delta W_i^n] = \delta_{ij}\tau$ and $\mathbb{E}[\Delta W_j^n\Delta W_i^{n+1}] = 0$ (because increments are independent). This yields
\begin{equation}
    \label{eq::LM-1}
    \begin{split}
    & g[\vec{Y}^n_{LM}, \tau_{n+1}]-g[\vec{Y}^n_{LM},\tau_n]= \partial_t g[\vec{Y}^n_{LM},\tau_n]\tau \\
    & +
    \nabla_{\vec x} g[\vec{Y}^n_{LM},\tau_n]^TA[\vec{Y}^n_{LM},\tau_n]\tau +  
    \frac{1}{2}\nabla g[\vec{Y}^n_{LM},\tau_n]^T B[\vec{Y}^n_{LM},\tau_n] (\Delta \vec{W}^n+\Delta \vec{W}^{n+1})+\\
    & + \frac{1}{4}\text{Tr}\{B B^T[\vec{Y}^n_{LM}(\tau_{n}),\tau_n] \text{Hess}_{\vec{x}}g[\vec{Y}^n_{LM}(\tau_{n}),\tau_n]\}\tau.   
    \end{split}
\end{equation}
Finally, using \eqref{Kolmogoroveq} for the partial derivative with respect to time we obtain (\ref{sol_BSDE_discret_LM}).

\section{Further assessment of the LM scheme for the heat equation}
\label{Appendix-B}
The Leihmkuhler-Matthews (LM) scheme appears to perform worse than both the Euler-Maruyama and Milstein schemes for the Black-Scholes equation (section \ref{sec:BS}). This is also the case when the diffusion tensor is homogeneous as tested on the Allen-Cahn equation (section \ref{sec:AC}). However we observe slight improvements when the time step $\tau$ of the discretization becomes smaller, and we therefore conclude that the LM scheme is limited by the number of neural networks needed for each time-slice. Taking smaller time steps adds networks which becomes prohibitive in terms of memory requirements.

However, the LM scheme is well suited for sampling from stationary distributions, and it is therefore natural to ask if this quality persists when solving an equation admitting a stationary limiting distribution. The simplest such equation is the heat equation that we discuss here. We find that for the same discretization interval used in the non-linear equations the LM scheme still leads to higher errors than Euler-Maruyama. 

\begin{figure}
    \centering
    \includegraphics[width=0.5\columnwidth]{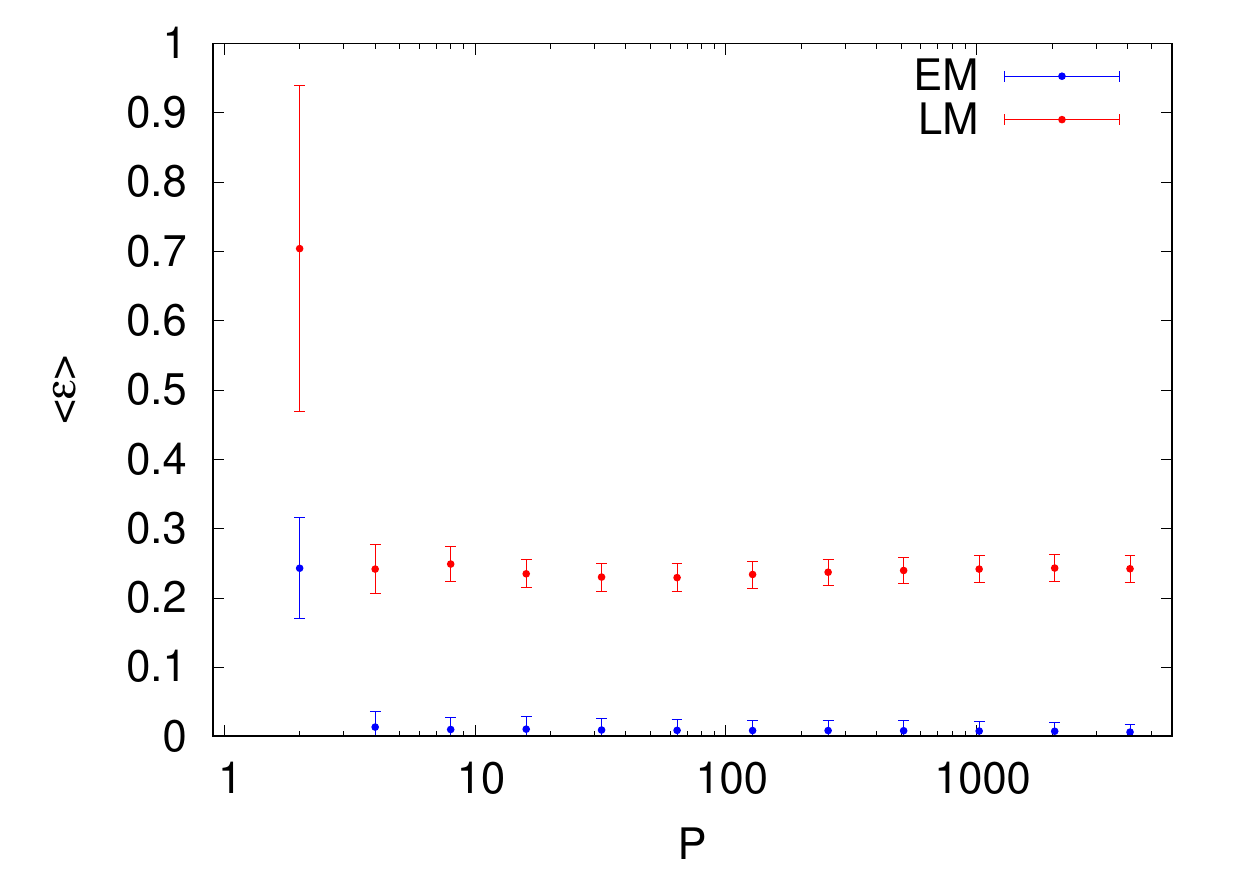}
    \caption{\small Comparison between Leimkuhler-Matthews (red points) and Euler-Maruyama (blue points) schemes for the simple {\it heat equation} in ten dimensions. The comparison is performed by computing the averaged relative approximation error defined in Eq. (\ref{rel_err_eq}) as a function of the number of initial points ($P$) for stochastic trajectoties used as inputs to the $\mathcal{DNN}$'s.  We observe that the Euler-Maruyama scheme is superior.}
    \label{EMvsLMrel_error}
\end{figure}

We look at the (forward) heat equation 
  $\frac{\partial g(\vec{x}, t)}{\partial t}=\Delta_x g(\vec{x}, t)$ 
in $d=10$ dimensions, on the interval $t \in [0,T]$, and with initial condition $g(\vec{x}, 0)=||\vec{x}||^2_{\mathbb{R}^d}$. The exact solution is given by $g(\vec{x}, t)=||\vec{x}||^2_{\mathbb{R}^d} + td$.

We only compare the  Leimkuhler-Matthews scheme with the Euler-Maruyama one (since the diffusion coefficient is constant). Fig. \ref{EMvsLMrel_error} shows the comparison between the two schemes by computing the averaged relative approximation error $\langle \epsilon \rangle$:
\begin{equation}\label{rel_err_eq}
    \langle \epsilon \rangle = \frac{1}{N} \sum_{i=1}^N \left| \frac{g(\vec{x}_i, T) - \mathcal{N}\mathcal{N}(\vec{x}_i, T|\vec{\theta})}{g(\vec{x}_i, T)} \right| 
\end{equation}
over $N=10^4$ points $\vec{x}_i$, $i=1,\dots,N$, randomly chosen in $[0,1]^d$, with $\mathcal{N}\mathcal{N}(\vec{x}, T|\vec{\theta})$ the solution returned by the respective $\mathcal{DNN}$. Fig. \ref{EMvsLMrel_error} shows that the Euler-Maruyama scheme is able to approximate better the solution of the heat equation than the Leimkuhler-Matthews one. 

\section{A simpler deep learning algorithm for linear equations. The example of geometric Brownian motion}
\label{Appendix-C}
The main focus of this paper is on non-linear Kolmogorov equations, however we briefly discuss for completeness a simpler algorithm that applies to linear equations. 
Physical examples for applications of linear Kolmogorv equations can be found in \cite{marino2016advective, marino2016entropy, aurell2016diffusion}.
For {\it linear} equations the standard Feynman-Kac formula is simply given by 
\begin{equation}
g(\vec{x}, t) =\mathbb{E}_{\vec x} [\phi(\vec{X}(T))]
\end{equation}
i.e., the first term in \eqref{eq:non-lin-FK}. In 
\cite{beck2018solving} it is beautifully remarked that this expectation minimizing a certain "mean square error" and that this can be taken as the basis of a deep learning algorithm. The main idea is to minimize the following loss:
\begin{align}\label{grosloss}
    \mathcal{L}(\vec{\theta})=\frac{1}{\vert b-a\vert^d}\int_{[a,b]^d} d\vec{x}\mathbb{E}_{\vec x}[
    \vert \phi(\vec{X}(T)) - \mathcal{N}\mathcal{N}(\vec{x}, T\vert \vec{\theta})\vert^2\bigr]
\end{align}
In practice the expectations are replaced by {\it empirical averages} over a sample set of discretized trajectories and the minimization over $\vec{\theta}$ carried on by a gradient descent algorithm. One important simplification with respect to the case of non-linear equations is that  only {\it one} deep network at time $T$ is optimized. 
We refer to \cite{beck2018solving} for further details of the algorithm. Our purpose here is to briefly compare the errors incurred by Euler-Maruyama and Milstein discretizations of the trajectories.

We work with a geometric Brownian motion, i.e., the linear part of the (backward) Black-Scholes equation \eqref{BS-eq}, 
\begin{equation}\label{eqGMB}
       \frac{\partial g(\vec{x},t)}{\partial t} + \frac{1}{2}\sum^d_{i=1}|\sigma_i x_i|^2 \frac{\partial^2 g}{\partial x^2_i}(\vec{x},t)+\sum^d_{i=1}\mu x_i\frac{\partial g}{\partial x_i}(\vec{x},t) = 0, \qquad g(\vec{x}, T) = \psi(\vec{x})
\end{equation}
for which the exact solution $g(\vec{x}, t)$, $0\leq t\leq T$ is known,
\begin{equation}
    g(\vec{x}, t)=\mathbb{E}[\psi\bigl(x_1 \exp(\sigma_1 W^{T-t}_1 + (\mu-\frac{|\sigma_1|^2}{2})(T-t)), \cdots, x_d \exp(\sigma_d W^{T-t}_1 + (\mu-\frac{|\sigma_1|^2}{2})(T-t))\bigr)]
\end{equation}
where $\vec{W}^t$ is the standard Brownian motion (conditioned to start at the origin $\vec{W}^0 = 0$). This is a one-dimensional integral that can be computed by the Monte Carlo method. 

The performance by computing the average relative error defined as 

\begin{equation}\label{rel_err_eq_int}
    \epsilon = \int_{[0,1]^d} \left| \frac{g(\vec{x}, 0) - \mathcal{NN}(\vec{x}, 0|\vec{\theta})}{g(\vec{x}, 0)} \right| d\vec{x}
\end{equation}
where $\mathcal{NN}(\vec{x}, 0|\vec{\theta})$, at $t=0$, is the approximate solution obtained by the deep learning algorithm. We take 
$r=\frac{1}{20}$, $\mu=r-\frac{1}{10}$, $\sigma_i=\frac{1}{10}+\frac{i}{200}$, $d=100$, $t\in [0,T=1]$, by using a time size interval $\tau=T/N$, with $N=40$. 
We fix the the initial condition to $g(\vec{x},0)=\psi(\vec{x})=\exp{(-rT)}\max\{[\max_{i \in {1,2,..,d}}x_i]-100,0\}$ (as in \cite{beck2018solving}).

\begin{figure}
    \centering
    \includegraphics[width=0.5\columnwidth]{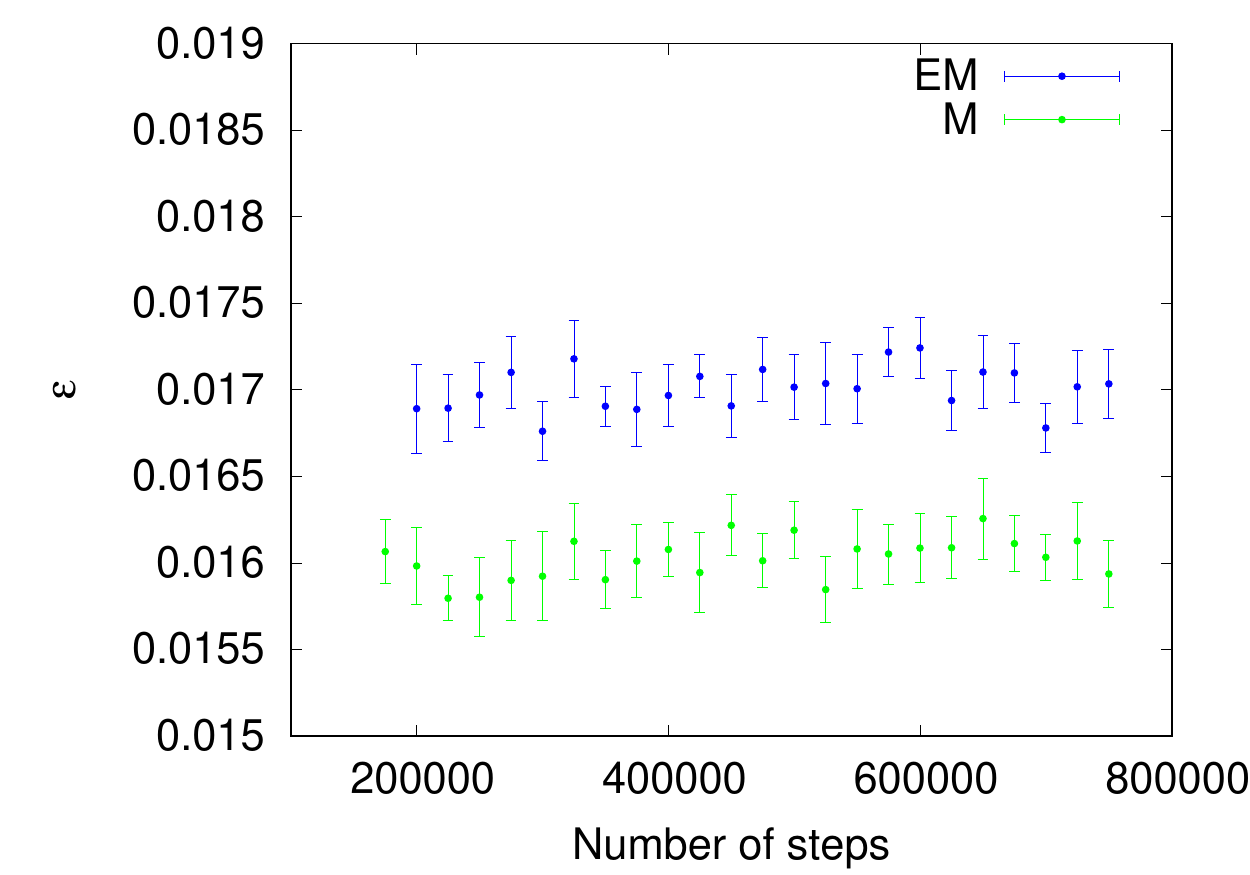}
    \caption{\small Comparison between Milstein (green points) and Euler-Maruyama (blue points) schemes. The comparison is performed by computing the relative approximation error defined in Eq. (\ref{rel_err_eq_int}) as a function of the number of steps done for training the parameters $\vec{\theta}$. The figure shows that the Milstein scheme better, on average, approximates the solution of the linear Black-Scholes equation as compared to the Euler-Maruyama one. Error bars are SEM.}
    \label{EMvsMrel_error}
\end{figure}

In Tab. $2$, and Fig. \ref{EMvsMrel_error} we report the results obtained for computing the relative error $\epsilon$ by using the two discretization schemes schemes in order to minimize the empirical averages corresponding to \eqref{grosloss}.  We also average $\epsilon$ over $15$ different experiments to reduce the variance. 
The analysis shows that exists an improvement in choosing an higher order discretization scheme for training the neural network. 

\begin{table}
\label{tab-d=8-9-10}
\centering
\scalebox{1}{
\begin{tabular}{ |c||c||c||} 
\hline
 Number of steps & $\epsilon_{EM} \pm \sigma_{mean}^{\epsilon_{EM}}$ &  $\epsilon_{M} \pm \sigma_{mean}^{\epsilon_{M}}$ \\
\hline
0 & 1.0001 (13)& 1.0018 (11) \\
25000 & 0.8427 (2)& 0.8421 (5)\\
50000 & 0.6855 (4) & 0.6847 (4)\\
75000 & 0.5290 (4) & 0.5282 (4)\\
150000 & 0.0804 (4) & 0.0777 (6)\\
200000 & 0.0169 (2) & 0.0159 (2)\\
500000 & 0.0170 (2) & 0.0162 (2)\\
750000 & 0.0170 (2) & 0.0159 (2)\\
\hline
\end{tabular}
}
\caption{\small The table shows the average values of relative error $\epsilon$. The integers in the brackets identify the standard errors of the mean on the last digits. On average Milstein approximation works better than the Euler-Maruyama one.  }
\end{table}


\section{Computational Complexity for the Milstein scheme}
\label{Appendix_complexity_Milstein}

In Sec. \ref{sec::non-linear-diff-eq}, we presented using the Euler-Maruyama scheme that the deep learning methodology studied in this paper does not suffer from the curse of dimensionality, in the sense that the number of training trajectories used scales polynomially in the dimension of space. In this appendix, however, we present the same analysis for the Milstein scheme, and we show that the number of training trajectories scales polynomially in the dimension of space. 

As far as we know, it is not know in literature an exact solution for a Kolmogorov non-linear equation with non constant transport coefficients. For this reason, we perform an accurate sampling of solutions using the multilevel Picard method for the equation (\ref{BSNL-eq}) in the interval $[49.995:50.005]^d$ and then we use those solutions as reference values for our analysis.

The deep network of Sec.\ref{Milstein_DNN} (with Milstein discretization) is used since the diffusion term is non-homogeneous. For the number of initial random points in $[49.995:50.005]^d$ of stochastic input trajectories we take $P=4,8,16,32, 64, 128, 256, 512, 1024$. Once the network is trained we obtain the 
approximate solution $\mathcal{N}\mathcal{N}(\vec{x},0|\vec{\theta})$.
The relative error is then easily computed from \eqref{rel_err_eq_DE}, using $M=10^2$ uniformly random points $\vec{x}_m\in [49.995:50.005]^d$, $m=1, \cdots, M$, computed with the multilevel Picard method. 

\begin{figure}
    \centering
    \includegraphics[width=0.45\columnwidth]{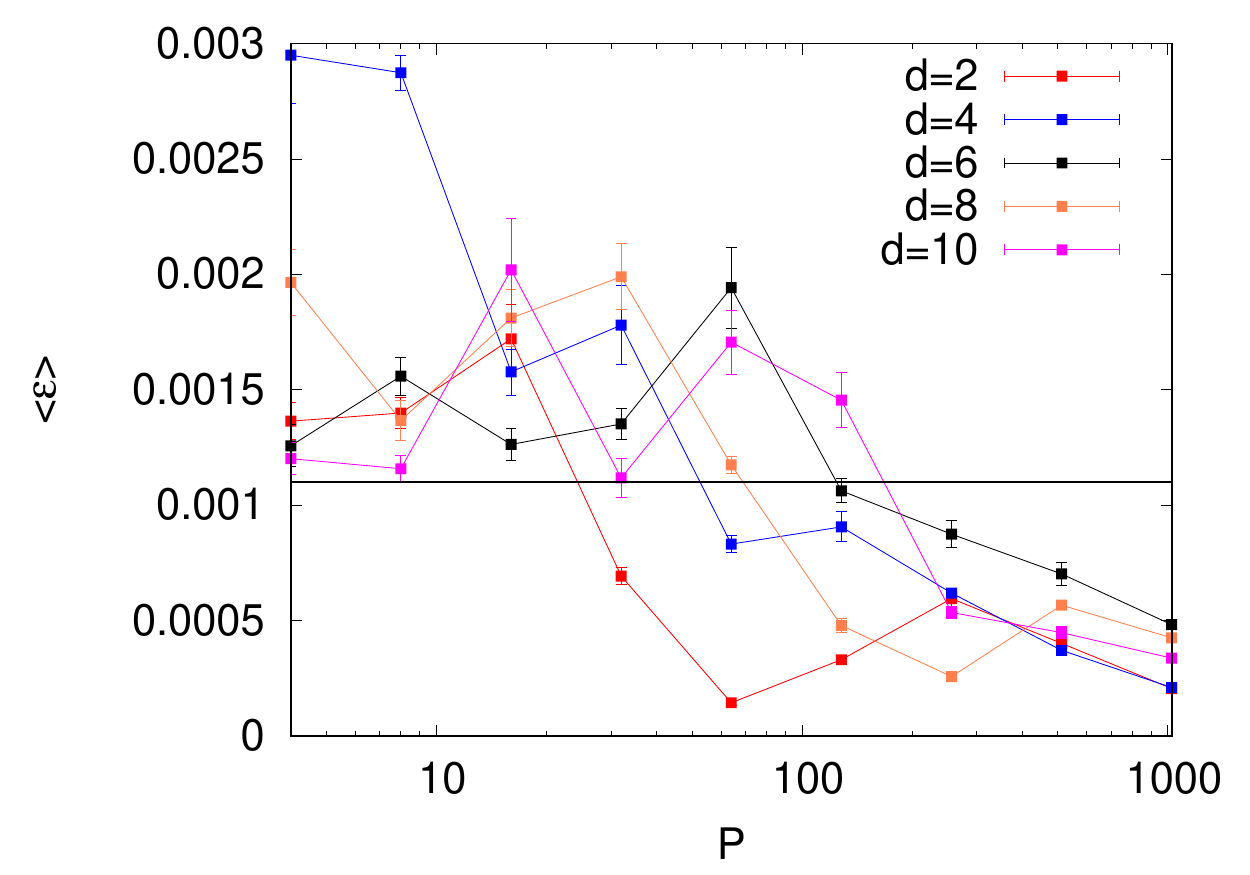}
    \includegraphics[width=0.45\columnwidth]{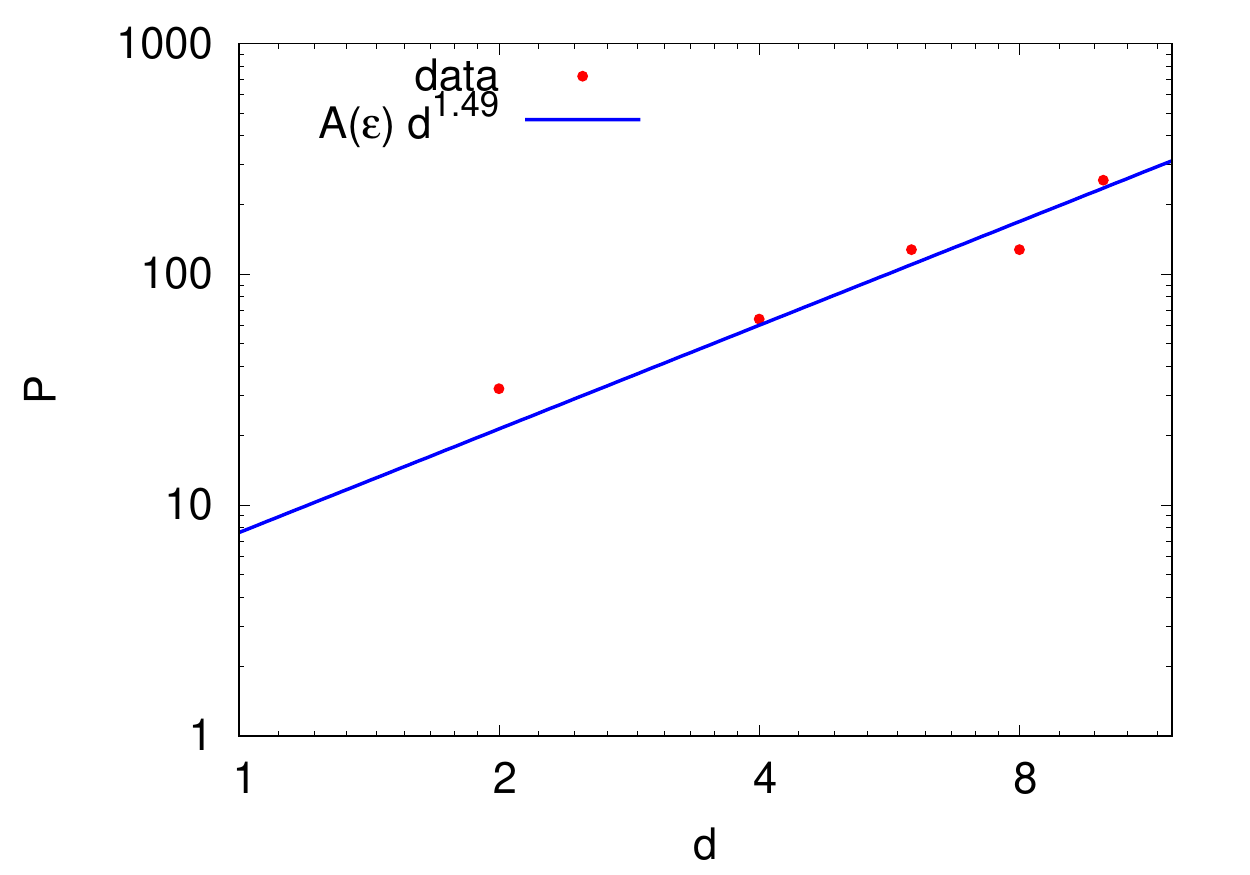}
    \caption{\small \textbf{Left}: Average relative error $\langle \epsilon \rangle$ as a function of number of initial points of stochastic trajectories $P=4,8,16,32, 64, 128, 256, 512, 1024$,
    on a log-log scale i.e. for different values of $d=2,4,6,8,10$.
    \textbf{Right}: Computational complexity, given an average relative error $\langle \epsilon \rangle < 0.0011 $, measured by $P$ as a function of $d=2,4,6,8,10$, on a log-log scale the slope is $b\approx 1.49$. The computational complexity obtained is proportional to $A(\epsilon)d^{1.49} $, with $A(\epsilon) \sim 7.6$. The results are obtained using the algorithm described in Sec.\ref{Milstein_DNN}, for the equation (\ref{BSNL-eq})  with $T=0.01$.}
    \label{fig-err-BSNL}
\end{figure}

Fig. \ref{fig-err-BSNL} shows the behaviour of $\langle\epsilon\rangle$ as a function of $P$, for various values of $d$. For small value of $d$, the averaged relative error needs few points $P$ to approximate well the exact solution. In contrast, as expected, as $d$ grows, the number of points $P$ needed for approximate well the solution becomes larger. 

By fixing a pre-specified relative error of $\langle \epsilon \rangle \approx 0.0011$  we observe that $P$ scales polynomially with $d$, roughly as $P=O(d^{1.49})$. This value of $P$ is comparable with the one obtained for the Euler-Maruyama scheme in Sec. \ref{sec::non-linear-diff-eq}. We conclude that the Milstein scheme improve the accuracy of the Neural Network, without modifying the computational complexity, as expected.  

\section{HJB equation with non constant diffusion coefficient}
\label{Appendix_HJBNCD}

In Sec. \ref{sec::HJB-sec}, we presented a comparison between Euler-Maruyama and Leihmkuhler-Matthews scheme over the Hamilton-Jacobi-Bellman equation in (\ref{HJB-eq}). In this appendix we present a comparison between the  Euler-Maruyama and a Mistein scheme over an Hamilton-Jacobi-Bellman equation with non constant diffusion tensor. The equation that we try to solve is:

\begin{equation}
\label{HJBNDC-eq}
    \frac{\partial g}{\partial t}(\vec{x},t) + \mathbf{D}(\vec{x})\Delta g(\vec{x},t)  = \lambda ||\nabla g(\vec{x},t)||^2,
\end{equation}
where $\lambda$ is a positive constant set to $1$, and $D(\vec{x})_{ij}=x_i^2 \delta_{ij}$, with $\delta_{ij}$ is the $\delta$-Kronecker. The terminal condition is $g(\vec{x}, T)=\phi(\vec{x})=\ln((1+||\vec{x}||^2)/2)$ with $\vec{x} \in \mathbb{R}^d$.

\begin{figure}
    \centering
    \includegraphics[width=0.9\columnwidth]{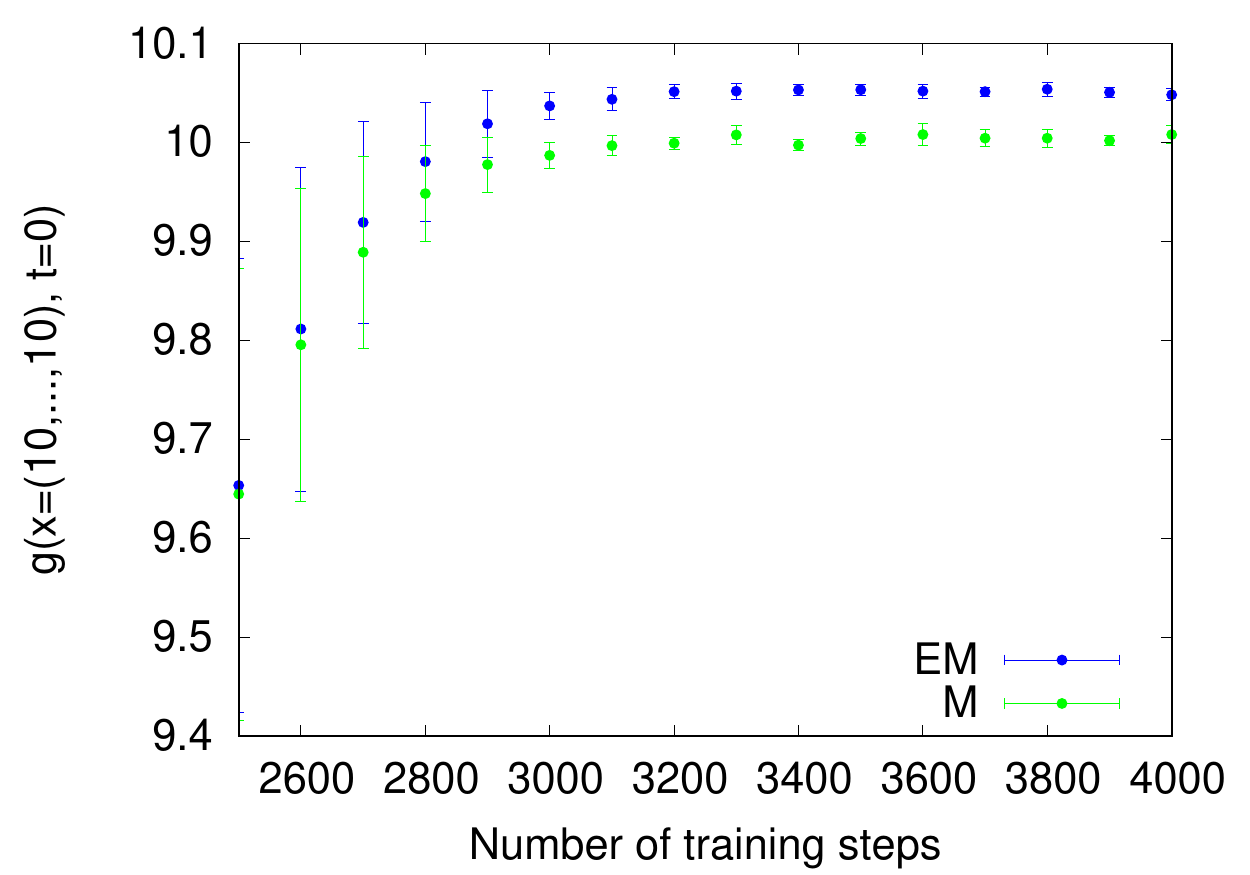}
    \caption{\small The figure shows the value of $g(\vec{x}=(50, \dots, 50), 0)$ returned by the neural networks $\mathcal{DNN}_{EM}$ (blue points), $\mathcal{DNN}_{M}$ (green points) as function of the number of training steps done for learning the parameters. Each point is the mean over five independent runs. The number of equidistant time steps was fixed at $40$, i.e. for each neural networks $\mathcal{DNN}_{EM}$ and $\mathcal{DNN}_{M}$  the value of $N$ is fixed to $40$. All the parameters were initialized randomly uniformly between $[-1,1]$. The total number of training steps, i.e. $t_s$, was fixed to $6000$ and the learning rate  to $\eta=0.008$.}
    \label{fig:HJBNL-comp}
\end{figure}

In Fig. \ref{fig:HJBNL-comp} we present the comparison between the output of $\mathcal{D}\mathcal{N}\mathcal{N}_{\rm EM}$ (blue points) and $\mathcal{D}\mathcal{N}\mathcal{N}_{\rm M}$ (green points). 
Each point is the average over 5 samples. Blue points identify the output of $\mathcal{D}\mathcal{N}\mathcal{N}_{\rm EM}$ and reach the approximate value $g_{\rm EM}(\vec{x}=(50, \dots, 50), 0)=10.049\pm 0.007$. The green points represent the value obtained by $\mathcal{D}\mathcal{N}\mathcal{N}_{\rm M}$ $g_{M}(\vec{x}=(50, \dots, 50), 0)=10.008 \pm 0.008$. Again, it is well shown the different results of the two schemes. No numerical or analytical result is present in the literature for equation (\ref{HJBNDC-eq}), but numerical observations show that the true value of the solution of $g(\vec{x}=(50, \dots, 50), 0)$ could be around $\sim 10.008$.

\section{Time scaling of the deep learning  algorithm}\label{Appendix-time}

In this appendix we briefly discuss how, in practice, the time scales with the dimension $d$, when the deep network is trained to compute an approximate solution. This experiment is performed on a cluster composed of 22 nodes, 200 cores and 250 GB of RAM divided in 22 different machines with different hardware.

We illustrate this for the heat equation and for the exactly solvable non-linear diffusion equation of Sec. \ref{sec::non-linear-diff-eq}. Since the diffusion tensor is constant we use the algorithm presented in Sec.\ref{subsec::EM} based on the Euler-Maruyama discretization. 

\begin{figure}[h!]
    \centering
    \includegraphics[width=0.40\columnwidth]{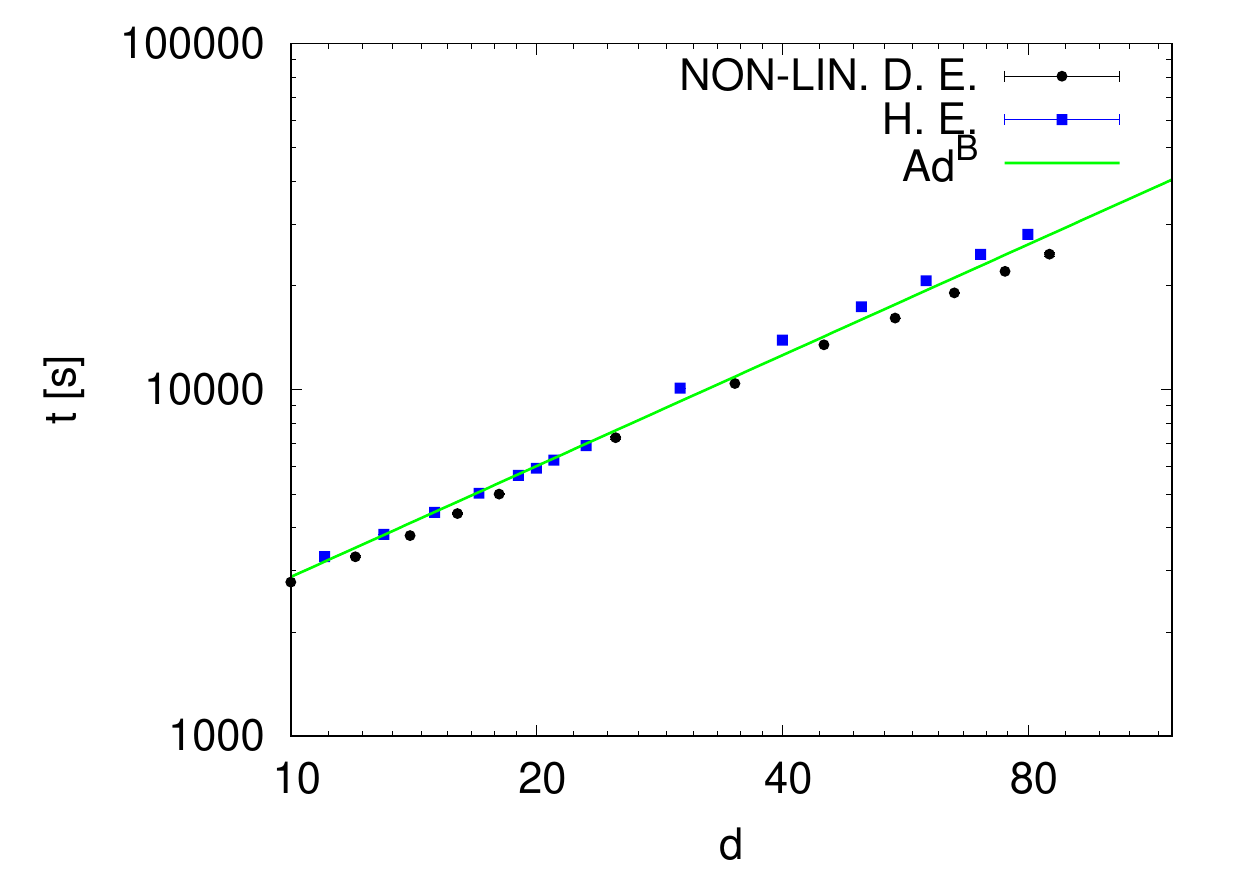}
    \caption{\small The figure shows how the time returned by the algorithm scales with the dimension $d$. Each point is an average over 3 experiments. Black circles correspond to the non-linear diffusion equation and blue squares to the heat equation. For both equations points live on the same straight line on a  log-log scale with slope $\sim 1.1$}. 
    \label{fig-time_cost}
\end{figure}

We use $P=4096$ random initial points for the set of stochastic trajectories used for training, and keep the number of training steps at $15000$. Fig. \ref{fig-time_cost} shows how the time needed to train the network and find an approximate a solution scales with the number of dimensions $d$. Interestingly we see that for both equations the set of points  live on the same straight line, on a log-log scale, with a slope $\sim 1.1$.

\ifCLASSOPTIONcaptionsoff
  \newpage
\fi



\bibliographystyle{IEEEtran}
\bibliography{IEEEabrv,references}
%
%
%
\end{document}